\documentclass{article}
\pdfoutput=1
\usepackage{amsmath}
\usepackage{amsfonts}
\usepackage{amssymb}
\usepackage{amsthm}
\usepackage[pdftex]{graphicx}
\usepackage{epstopdf}
\usepackage{cite}
\usepackage{makecell}
\usepackage{bm}%
\providecommand{\U}[1]{\protect\rule{.1in}{.1in}}
\providecommand{\U}[1]{\protect\rule{.1in}{.1in}}
\usepackage[margin=0.8in]{geometry}
\newtheorem{theorem}{Theorem}

\newtheorem{definition}[theorem]{Definition}

\newtheorem{proposition}[theorem]{Proposition}
\newtheorem{remark}[theorem]{Remark}

\begin{document}
	
	\title{\LARGE \textbf{Exact recursive updating of uncertainty sets}}
	\author{Robin Hill$^{1}$, Yousong Luo$^{2}$ and Uwe Schwerdtfeger$^{3}$ \thanks{$^{1}%
			$Department of Electrical and Electronic Engineering, University of Melbourne,
			Melbourne, Vic 3010, Australia \texttt{{\small robin.hill@unimelb.edu.au}}%
		}\thanks{$^{2}$School of Science, RMIT
		University, 124 Latrobe St, Melbourne, 3001, Australia
		\texttt{{\small yluo@rmit.edu.au}}}\thanks{$^{3}$Department of Mathematics,
		Chemnitz University, Germany
		\texttt{{\small uwe.schwerdtfeger@mathematik.tu-chemnitz.de}}}}
\maketitle

\begin{abstract}                          
This paper addresses the classical problem of determining the set of possible
states of a linear discrete-time system subject to bounded disturbances from
measurements corrupted by bounded noise. These so-called uncertainty sets
evolve with time as new measurements become available. We present two theorems which describe completely how they evolve with time, and this yields an efficient algorithm for recursively updating uncertainty sets. Numerical simulations demonstrate performance improvements over existing exact methods.  

\end{abstract}

\section{INTRODUCTION}
\setcounter{page}{1}
Consider a linear, time-invariant dynamic system driven by set-bounded process
noise, and with measurements corrupted by set-bounded observation noise. The set of possible states of the system consistent with the
measurements up to the current time is termed the \textit{state uncertainty
	set }(or simply \textit{uncertainty set}). In many applications having a representation of the uncertainty set is useful. This so-called set membership estimation problem
is fundamental and has many applications, for example in fault detection \cite{tabatab2015,Rosa-CDC2010,TORNILSIN20121,alamo-etal-2005,Casau_etal-2015}, control under
constraints in the presence of noise \cite{Bertsekas-etal-1971,glover-schweppe-1971}, 
and 
model (in)validation \cite{Poolla_etal-1994,Rosa_silvestre_athans-2014}. A closely related topic is identification of bounded-parameter models \cite{Belforte_etal_1990,Clement_Gentil_1990,Norton_1987}. 

The first results on recursive determination of the uncertainty set are in
\cite{Schweppe-1968} and
\cite{Witsenhausen-poss-states-1968}.
Since the appearance of these papers there has appeared an extensive literature on
the topic. See \cite{Fogel-Huang-82} and \cite{Ninness_Goodwin-95} for background on the set-bounded approach to uncertainty, the survey paper \cite{Milanese-Vicino-91} and the book \cite{blanchini_miani}. Some of the many other papers which consider this problem are \cite{Blanchini-Sznaier-2012,Stoorvogel-1996,Tempo-1988}.

In the first part of the seminal paper \cite{Witsenhausen-poss-states-1968} an exact in principle solution to the problem of recursively determining polytopic uncertainty sets is given. It uses the $\mathcal{H}$-representation for the uncertainty sets,  that is they are defined using inequality constraints. But the solution requires (Minkowski) addition, and intersections, of polytopes, both of which can be time-consuming. Exact, recursive $\mathcal{H}$-representation methods often use Fourier-Motzkin elimination or parametric linear programming, see \cite{Keerthi_Gilbert-1987,Rakoviv-Mayne-2004,Shamma_Kuang-1999} for the former, and \cite{Jones2008} for the latter. In these implementations it is the identification and removal of redundant inequality constraints that is most demanding computationally. The redundant constraints can be removed by solving linear programs but this is not a trivial task, for which only weak polynomiality is known if only the $\mathcal{H}$-representation of the polytope is available. For hardness results on polytopic computations, see \cite{Tiwary2008}.

Another interesting recent approach using exact methods, based on geometric ideas, is in \cite{hagemann}. Here also an inequality description is used, and projection followed by redundant inequality constraint elimination is necessary.  In our algorithm we are in the fortunate position of having both vertex and inequality representations. This means we can efficiently intersect hyperplanes and polytopes, as well as pairs of facets, the only computationally intensive tasks our algorithm requires.

In Section V of \cite{Witsenhausen-poss-states-1968} a dual to the $\mathcal{H}$-representation is presented. Using the theory of conjugate functions an equation describing dynamic evolution of the support function of the uncertainty set is derived, see Section VII, page 558, for the special case of independently constrained noise signals, the case considered in this paper. While of very significant theoretical value, the results in \cite{Witsenhausen-poss-states-1968} were not developed to the point of yielding an algorithm for uncertainty set propagation. 

In this paper we build on the ideas in \cite{Witsenhausen-poss-states-1968}, particularly that of support function evolution. We use linear programming rather than conjugate functions as our basic tool, and employ the familiar complementary slackness conditions relating primal and dual variables to prove our main results. 

In real-time applications, for example fault detection and isolation \cite{Rosa-CDC2010}, existing exact algorithms run up against the problem of their computational complexity. For this reason there has been a lot of research recently on the use of zonotopes and constrained zonotopes to approximate the exact polytopic uncertainty set, see for example \cite{Combastel2015265,Alamo-etal-2008,Scott2016126}. 

The results in this paper provide tools for investigating how the complexity of the polytopic uncertainty set varies as more measurements become available.
For example, the complexity of the polytopic representation of the uncertainty sets for a fifth order plant is very variable, but does not appear to have long term growth when the measurements are randomly selected. For higher order plants the growth in complexity is faster and it is not yet clear if there is ever any levelling off in complexity.

If approximation is necessary, our exact results could also be useful. Having an exact representation enables intelligent approximation. For example it has been noticed in our simulations that vertices often accumulate close to each other, on facets having almost identical directions. Identifying such behaviour allows for greater complexity reduction with smaller error. There is perhaps scope for combining exact and zonotope approximations in the trade off between complexity and efficiency.


\section{Basic Setup}

The plant $P$, a linear, time-invariant, causal discrete-time, $m^{\rm th}$
order scalar system, is assumed known. There are two sources of
uncertainty, an input noise disturbance $(u_{k})_{k=1}^{\infty}=\mathbf{u},$
and output measurement noise $(w_{k})_{k=1}^{\infty}=\mathbf{w}.$ The plant
output is $(y_{k})_{k=1}^{\infty}=\mathbf{y}$, and the measurement at time $k$
is $z_{k}=y_{k}+w_{k}.$ The initial state, at time $k=0,$ is assumed to be
known exactly, but nothing is known about the uncertainties except that they
satisfy $\left\vert u_{k}\right\vert \leq1$ and $\left\vert w_{k}\right\vert
\leq1.$ We will refer to this as the primal system.

Given an initial state $\mathbf{x}_{0},$ the measurement history $z_{1}%
,\ldots,z_{k-1},$ and the plant dynamics, we seek the uncertainty set
at time $k,$ denoted $S_{k};$ it is the set of possible states at time $k$ consistent with
the measurements up to and including $z_{k-1}$, and is easily seen to be a closed, convex polytope.

\subsection{Notation\label{sectnotprelim}}

Given a vector $\mathbf{y=}\left(  y_{1},y_{2},\ldots\right)  $ and any
$s\in\mathbb{N}^{+},$ $t\in\mathbb{N}^{+}$ satisfying $s<t,$ we denote
$\left(  y_{s},y_{s+1},\ldots,y_{t}\right)  $ by $y_{s:t}.$ The $\lambda
$-transform (generating function) of an arbitrary sequence $\mathbf{y}%
=(y_{k})_{k=1}^{\infty}$ is defined to be $\hat{\mathbf{y}}(\lambda):=\sum
_{k=1}^{\infty}y_{k}\lambda^{k-1}.$  Real Euclidean space of dimension $m$ is denoted $\mathbb{R}^{m}$, where $m$ is the order (McMillan degree) of the plant $P$. States of the plant $P$ are represented by vectors, or points, in  $\mathbb{R}^{m}$.  Let $\mathbf{d}=d_{1:m+1}=(d_{1}%
,\ldots,d_{m+1})$ and $\mathbf{n}=n_{1:m+1}=(n_{1},\ldots,n_{m+1}),$ 
be real vectors, where $\hat{\mathbf{n}}(\lambda)$ and $\hat{\mathbf{d}}(\lambda)$ are the numerator and denominator of the transfer function representation of the plant $P$. Denote by $\mathbf{D}_{\infty}$ and $\mathbf{N}_{\infty}$ the infinite, banded, lower-triangular
Toeplitz matrices whose first columns are $\mathbf{d}$ and $\mathbf{n}$,
respectively. Define the following lower and upper triangular submatrices of $\mathbf{D}_{\infty}.$%
\begin{align*}
\mathbf{D}_\mathrm{L}  &  :=\left[
\begin{array}
[c]{cccc}%
d_{1} & 0 & \dotsc & 0\\
d_{2} & d_{1} & \ddots & \vdots\\
\vdots & \ddots & \ddots & 0\\
d_{m} & \dotsc & d_{2} & d_{1}%
\end{array}
\right] \\
\mathbf{D}_\mathrm{U}  &  :=\left[
\begin{array}
[c]{cccc}%
d_{m+1} & d_{m} & \dotsc & d_{2}\\
0 & d_{m+1} & \ddots & \vdots\\
\vdots & \ddots & \ddots & d_{m}\\
0 & \dotsc & 0 & d_{m+1}%
\end{array}
\right]  .
\end{align*}
The matrices $\mathbf{N}_\mathrm{L}$ and $\mathbf{N}_\mathrm{U}$ are defined similarly.

For any $k>0,$ the $k\times k$ upper left hand corner submatrix of $\mathbf{D}_{\infty}$ is
denoted $\mathbf{D}_{k \times k}.$ We will often write simply $\mathbf{D}$ instead of $\mathbf{D}_{k \times k}$ when $k$ is clear from context. The symbols $\mathbf{N}_{k \times k}$ and $\mathbf{N}$ are defined similarly. Note that $\mathbf{D}_{m \times m}=\mathbf{D}_\mathrm{L}$
and $\mathbf{N}_{m \times m}=\mathbf{N}_\mathrm{L}.$

The Toeplitz Bezoutian matrix of $\mathbf{n}$ and
$\mathbf{d}$ is defined as $\mathbf{B}_\mathrm{T}:=\mathbf{D}_\mathrm{L}\mathbf{N}_\mathrm{U}-\mathbf{N}_\mathrm{L}\mathbf{D}_\mathrm{U}$.

One form of the Gohberg-Semencul formulas \cite{gohberg-semencul-1972,fuhrm} states%

\begin{equation}
\mathbf{B}_\mathrm{T}=\mathbf{N}_\mathrm{U}\mathbf{D}_\mathrm{L}-\mathbf{D}_\mathrm{U}\mathbf{N}_\mathrm{L},
\label{def:BTmatrix}%
\end{equation}
and this will be needed in the proof of Theorem \ref{main}, which underpins all of our results.
The first row of $\mathbf{B}_\mathrm{T}$ plays an important role and will be denoted
by $\mathbf{C},$ so 
$\mathbf{C}:= d_{1}\left[n_{m+1}, \ldots,n_{2}\right] -n_{1}\left[d_{m+1}, \ldots, d_{2}\right]$.

The inverse of $\mathbf{B}_\mathrm{T}$ exists if the polynomials $\hat{\mathbf{n}}(\lambda)$
and $\hat{\mathbf{d}}(\lambda)$ are
coprime, and $\mathbf{B}_\mathrm{T}^{-1}$ denotes the inverse of $\mathbf{B}_\mathrm{T}.$
See \cite{Heinig-Rost-2010} for properties of Bezoutians.

\subsection{Transfer function description and state-space representations}

The plant for the primal system has the transfer function representation
$P(\lambda)=\hat{\mathbf{n}}(\lambda)/\hat{\mathbf{d}}(\lambda)$ where
\begin{align}
\hat{\mathbf{n}}(\lambda)  &  =n_{1}+n_{2}\lambda+n_{3}\lambda^{2}%
+\dotsc+n_{m+1}\lambda^{m}\nonumber\\ 
\hat{\mathbf{d}}(\lambda)  &  =d_1+d_{2}\lambda+d_{3}\lambda^{2}+\dotsc
+d_{m+1}\lambda^{m},\nonumber
\end{align}
$m\geq1$ is an integer, $\hat{\mathbf{n}}(\lambda)$ and $\hat{\mathbf{d}%
}(\lambda)$ are assumed to be coprime polynomials with real coefficients, and
it is assumed that both the plant $P(\lambda)$ and the plant $P^{\ast}%
(\lambda)$ for the dual system, defined below, are causal, implying
$d_{1}\neq0$ and $d_{m+1}\neq0$, in which case the system matrix $\mathbf{A}$ is non-singular. 
Without loss of generality we take $d_{1}=1.$
Assuming zero initial conditions, $\mathbf{y}$ and $\mathbf{u}$ are related
by $\hat{\mathbf{d}}(\lambda)\hat{\mathbf{y}}(\lambda)=\hat{\mathbf{n}}(\lambda)\hat{\mathbf{u}}(\lambda).$ 

A state-space description of the primal system is
\begin{align}
\mathbf{x}_{k+1}  &  =\mathbf{A}\mathbf{x}_{k}+\mathbf{B}u_{k}\label{primalstate1}\\
y_{k}  &  =\mathbf{C}\mathbf{x}_{k}+Du_k\label{primalstate2}\\
z_{k}  &  =y_{k}+w_{k}\nonumber
\end{align}

where
\begin{align} 
\mathbf{A}  &  =\left[
\begin{array}
[c]{cc}%
\mathbf{0} & \mathbf{I}_{m-1}   \\
-d_{m+1} &
\begin{array}
[c]{cc}%
\dotsc & -d_{2}%
\end{array}
\end{array} 
\right]  ,\text{ }\mathbf{B}=\left[
\begin{array} 
[c]{c}%
\mathbf{0}  \\
1
\end{array}
\right]  , \label{ABCD} \\ \nonumber 
\mathbf{C}  &  = d_{1}\left[n_{m+1}, \ldots,n_{2}\right] -n_{1}\left[d_{m+1}, \ldots, d_{2}\right], \text{ }D=n_1, \nonumber  \\ \nonumber 
\end{align}
and the state at time $k \ge 0$ for the sequence pair $\left( \mathbf{y},\mathbf{u} \right) =\left( (y_j)_{j=1}^{\infty},(u_j)_{j=1}^{\infty}\right)$ is given by

\begin{equation}
\mathbf{x}_{k}:=\left\{
\begin{array}
[c]{c}%
\mathbf{B}_\mathrm{T}^{-1}\left[  \mathbf{D}_\mathrm{L}y_{k+1:k+m}-\mathbf{N}_\mathrm{L}u_{k+1:k+m}\right]  \text{ for }%
k\geq0\\
\mathbf{B}_\mathrm{T}^{-1}\left[  -\mathbf{D}_\mathrm{U}y_{k-m+1:k}+\mathbf{N}_\mathrm{U}u_{k-m+1:k}\right]  \text{ for }k\geq
m.
\end{array}
\right.  \label{statedef}%
\end{equation}

There is a system closely related to the estimation system that we refer to
as the dual system. The dual system input and output sequences are $(y_{j}^{\ast})_{j=1}^{\infty}$
and $(u_{j}^{\ast})_{j=1}^{\infty},$ and the dual plant, denoted
$P^{\ast},$ has the transfer function representation%

\begin{equation}
P^{\ast}(\lambda)=-\hat{\mathbf{n}}_{\mathrm{dual}} (\lambda)/ \hat{\mathbf{d}}_{\mathrm{dual}}(\lambda) 
\label{regconv}%
\end{equation}
where $\mathbf{n}_{\mathrm{dual}}=\left(  n_{m+1},\ldots,n_{1}\right)  $ and
$\mathbf{d}_{\mathrm{dual}}=\left(  d_{m+1},\ldots,d_2,1\right)$. A minimal
state-space realization of the dual system is%

\begin{align}
\mathbf{x}_{k+1}^{\ast}  &  =\mathbf{A}^{\ast}\mathbf{x}_{k}^{\ast}+\mathbf{B}^{\ast}y_{k}%
^{\ast}\label{adjointstate1}\\
u_{k}^{\ast}  &  =\mathbf{C}^{\ast}\mathbf{x}_{k}^{\ast}+D^{\ast}y_{k}^{\ast}
\label{adjointstate}%
\end{align}%
\begin{align}
\mathbf{A}^{\ast}  &  =\left[
\begin{array}
[c]{cc}%
-d_{m}/d_{m+1} & \mathbf{I}_{m-1}\\%
\begin{array}
[c]{c}%
\vdots\\
-1/d_{m+1}%
\end{array}
&
\begin{array}
[c]{c}%
\\
\mathbf{0}
\end{array}
\end{array}
\right]  ,\label{ABCDstardef}\\
\text{ }\mathbf{B}^{\ast}  &  =\left[
\begin{array}
[c]{c}%
n_{m}\\
\vdots\\
n_{1}%
\end{array}
\right]  -\left[
\begin{array}
[c]{c}%
d_{m}\\
\vdots\\
1
\end{array}
\right]  \frac{n_{m+1}}{d_{m+1}},\label{ABCDstardef1}\\
\mathbf{C}^{\ast}  &  =\left[
\begin{array}
[c]{cccc}%
-1/d_{m+1} & 0 & \dotsc & 0
\end{array}
\right]  ,\text{ }D^{\ast}=\frac{-n_{m+1}}{d_{m+1}};\label{ABCDstardef2}\\\nonumber
\end{align}
and the dual state at time $k \ge 0$ for the sequence pair $\left( \mathbf{y}^{\ast},\mathbf{u}^{\ast} \right) =\left( (y_j^{\ast})_{j=1}^{\infty},(u_j^{\ast})_{j=1}^{\infty}\right)$ is given by
\begin{equation}
\mathbf{x}_{k}^{\ast}:=\left\{
\begin{array}
[c]{c}%
-\mathbf{N}_\mathrm{U}^\mathrm{T}y_{k+1:k+m}^{\ast}-\mathbf{D}_\mathrm{U}^\mathrm{T}u_{k+1:k+m}^{\ast}\text{ for }k\geq0\\
\mathbf{N}_\mathrm{L}^\mathrm{T}y_{k-m+1:k}^{\ast}+\mathbf{D}_\mathrm{L}^\mathrm{T}u_{k-m+1:k}^{\ast}\text{ for }k\geq m
\end{array}
\right.  , \label{state_def_star}%
\end{equation}
The primal and dual state-space representations are
in principle well known \cite{kailath_linsys,LUEN-1979,pold_willems}. They are given explicitly in \cite{Luo_Hill-2015}.

\subsection{Polytopes}
The primal and dual states defined above will be interpreted in terms of the geometry of the polytopic uncertainty sets, so in this section we introduce notation and briefly summarise the relevant theory of convex polytopes. For more information and background, including the definition of a polytope, see for example 
\cite{Lay_convsets}, \cite{Acary-etal-2011} or \cite{ziegler_polytopes}.
Let $S$ denote a closed, convex polytope.

The \textit{support function} of 
$S$ in $\mathbb{R}^{m}$ is 
\[h_{S}(\mathbf{f}%
)=\max_{\mathbf{x}\in S}\left\langle \mathbf{f},\mathbf{x}\right\rangle
,\] where $\mathbf{f}\in\mathbb{R}^{m}.$

When $\mathbf{f} \ne \mathbf{0}$ the set
\[
H_{S}(\mathbf{f}):=\left\{  \mathbf{x}\in\mathbb{R}^{m}:\left\langle
\mathbf{f},\mathbf{x}\right\rangle =h_{S}(\mathbf{f})\right\}
\]
is the supporting hyperplane of $S$ with direction (outer normal vector)
$\mathbf{f}.$  If $\mathbf{f}=\mathbf{0}$ then $H_{S}(\mathbf{f})=\mathbb{R}^{m}$.

The intersection of $S$ with a supporting hyperplane is
called a \textit{face} of $S,$ and a face of dimension $m-1$ is a \textit{facet} of $S$.  A face of dimension $m-2$ is called a \textit{ridge}, and the faces of dimensions $0$ and $1$ are termed \textit{vertices} and \textit{edges}, respectively.

\begin{figure}[ptb]
	\centering
	\includegraphics[scale=0.45]{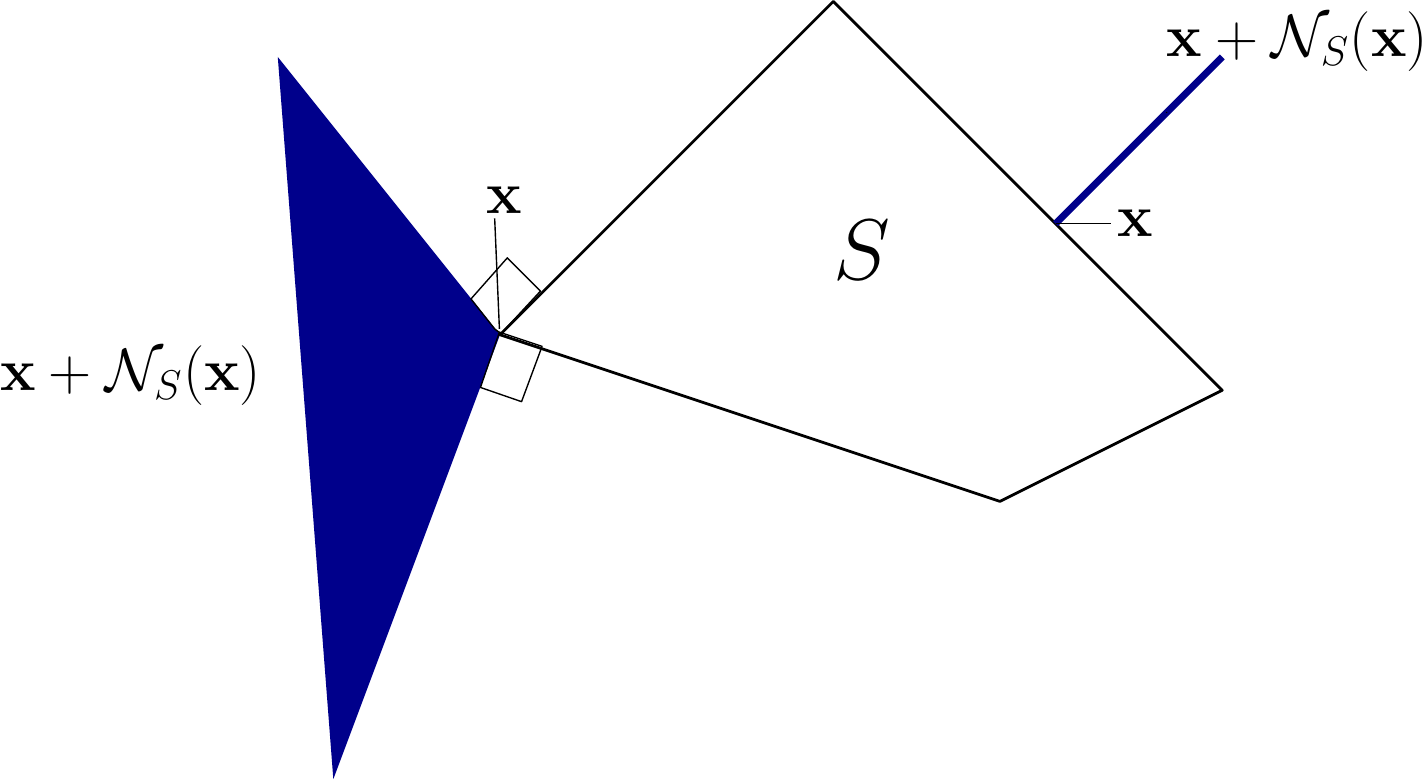} \caption{Normal cones to the polytope $S$}%
	\label{ncones}%
\end{figure}

The \textit{normal cone} at a boundary point $\mathbf{x}$ is the set \[\left\{
\mathbf{f}\in\mathbb{R}^{m}:\left\langle \mathbf{f},\mathbf{x}\right\rangle
=h_{S}(\mathbf{f})\right\}  \]
and is denoted $\mathcal{N}_S\left( \mathbf{x}\right)$.
It is generated by the outward normals to the facets that form the polytope at $\mathbf{x}$, that is
\[
\mathcal{N}_{S}\left( \mathbf{x}\right)=\{\alpha_1\mathbf{f}_1+\dotsc+\alpha_n \mathbf{f}_n:\alpha_1,\ldots,\alpha_n\geq0\},
\]
where $\mathbf{f}_1,\ldots,\mathbf{f}_n$ are the directions of the facets containing $\mathbf{x}$. 
Thus $\mathcal{N_S}\left(
\mathbf{x}\right)$ contains the directions of all
hyperplanes which touch $S$\ at $\mathbf{x.}$ See Fig. \ref{ncones}. If $\mathbf{x} \in \operatorname*{int}\left(S\right)$, then $\mathcal{N_S}\left(%
\mathbf{x}\right):=\left\{\mathbf{0}\right\}$. By definition the directions of facets, and of the hyperplanes that contain facets, point outwards from the polytope. A direction is a non-zero vector. 

In Section \ref{sectprimaldual} the dual state $\mathbf{x}^*_k\left(\mathbf{y}^*,\mathbf{u}^*\right)$ will be interpreted as a  vector $\mathbf{f}$ in the normal cone of the primal state $\mathbf{x}_k\left(\mathbf{y},\mathbf{u}\right) \in S_k$. The symbol $\mathbf{f}$ will be used when the geometric viewpoint is being emphasised, while $\mathbf{x}^*_k\left(\mathbf{y}^*,\mathbf{u}^*\right)$ will be used to denote the same vector from the system theoretic, algebraic point of view. We shall sometimes drop the subscript $k$, and indicate the next time instant with the superscript $^+$. For example, $\mathbf{f}^+$ will denote $\mathbf{x}^*_{k+1}\left(\mathbf{y}^*,\mathbf{u}^*\right)$.
The alternative notations are summarised below:
\begin{flushleft}
	\begin{equation}
	\begin{array}
	[c]{c}%
	\mathbf{f}\leftrightarrow\mathbf{x}^*_k\left(\mathbf{y}^*,\mathbf{u}^*\right)\\
	\mathbf{f}^+\leftrightarrow\mathbf{x}^*_{k+1}\left(\mathbf{y}^*,\mathbf{u}^*\right)\\
	\mathbf{x} \leftrightarrow \mathbf{x}_k \left(\mathbf{y},\mathbf{u}\right)\\
	\mathbf{x}^+ \leftrightarrow \mathbf{x}_{k+1} \left(\mathbf{y},\mathbf{u}\right).
	\end{array}
	\nonumber 
	\end{equation}
\end{flushleft}

\subsection{State propagation\label{sect_prop_uncset}}

The primal system at time zero is in the state $\mathbf{x}_{0}$ so, by (\ref{statedef}),
$\mathbf{D}_\mathrm{L}y_{1:m}-\mathbf{N}_\mathrm{L}u_{1:m}=\mathbf{B}_\mathrm{T}\mathbf{x}_{0}.$ At any time $k>m$, $y_{1:k}$ and $u_{1:k}$ are related by%
\begin{equation}
\mathbf{D}y_{1:k}-\mathbf{N}u_{1:k}=\left[
\begin{array}
[c]{c}%
\mathbf{B}_\mathrm{T}\mathbf{x}_{0} \\
\mathbf{0}
\end{array}
\right]. \label{conu_indconf}%
\end{equation}
Equation (\ref{conu_indconf}) is a consequence of the plant input/output relationship in transfer function form:
\begin{equation}
\hat{\mathbf{d}}\hat{\mathbf{y}}-\hat{\mathbf{n}}\hat{\mathbf{u}}=\hat{\mathbf{b}}, \label{primeconsts}
\end{equation}
where we have used the abbreviation $\mathbf{b}=\mathbf{B}_\mathrm{T}\mathbf{x}_{0}$. Then (\ref{conu_indconf}) follows from equating like powers of $\lambda$ on both sides of (\ref{primeconsts}).

By the state-space representation of the primal system,  $\mathbf{x}_{k+1}=\mathbf{A}\mathbf{x}_{k}+\mathbf{B}u_{k}$ and $y_{k}  =\mathbf{C}\mathbf{x}_{k}+Du_{k}$. Recall that $S_k$ is the set of states at time $k$, given measurements up to time $k-1$. Following Witsenhausen, \cite{Witsenhausen-poss-states-1968}, the set $S_{k+1}$ is given
recursively in terms of $S_{k}$ and the new observation $z_{k}$ by

\begin{equation} 
S_{k+1}=\left\{ \mathbf{x}_{k+1} \; : \;
\begin{array}
[c]{l}%
\mathbf{x}_{k}\in S_{k}, \mathbf{x}_{k+1}=\mathbf{A}\mathbf{x}_{k}+\mathbf{B}u_{k},\\
y_{k}=\mathbf{C}\mathbf{x}_{k}+Du_{k},\\
\left\vert u_{k}\right\vert \leq1,\text{ }\left\vert y_{k}-z_{k}\right\vert
\leq1.
\end{array} 
\right\}   \label{wits_recurs}%
\end{equation}

Special notation is now introduced to describe states $\mathbf{x}_{k+1}$ and $\mathbf{x}_{k}$ related as in (\ref{wits_recurs}).

\begin{definition}
	\label{defprecursor}The state $\mathbf{x}_{k} \in S_{k}$ is said to be a
	\textit{precursor} of the state $\mathbf{x}_{k+1}$, $\mathbf{x}_{k}$
	is \textit{propagated} to $\mathbf{x}_{k+1},$ and $\mathbf{x}_{k+1}$ is a \textit{successor} to
	$\mathbf{x}_{k},$ if there exists a scalar $u_k$ satisfying $\left\vert u_{k}\right\vert \leq1$ for which $\mathbf{x}_{k+1}=\mathbf{A}\mathbf{x}_{k}+\mathbf{B}u_{k}$, and $\left\vert y_{k}-z_{k}\right\vert
	\leq1$ where $y_{k}=\mathbf{C}\mathbf{x}_{k}+Du_{k}$.
\end{definition}
So $S_{k+1}$ is the set of all successors to all states in $S_{k}$, and any precursor of any state $\mathbf{x}_{k+1}\in S_{k+1}$ is in $S_{k}$.

Associated with any state $\mathbf{x} \in S_{k}$ define in the
$y_{k}u_{k}$-plane the \textit{primal line} $L\left(\mathbf{x}\right)$:
\begin{definition}
	$L\left(\mathbf{x}\right)=\{(y_k,u_k): y_{k}-n_{1}u_{k}=\mathbf{C}\mathbf{x}\}.$
\end{definition}
The input and output of the plant at time $k$ are constrained to lie on the line $L\left(  \mathbf{x}\right)$ by (\ref{primalstate2}) and (\ref{ABCD}). 
Associated with any measurement $z_k$ define a square $Q$ in the $y_ku_k$-plane.
\begin{definition}
\label{defQ}
	$
	Q(z_k)=\{(y_k,u_k):\left\vert u_{k}\right\vert \leq1 \text{  and  } \left\vert y_{k}-z_{k}\right\vert
	\leq1 \}.
	$
\end{definition}

\begin{figure}[ptb]
	\centering
	\includegraphics[scale=0.5]{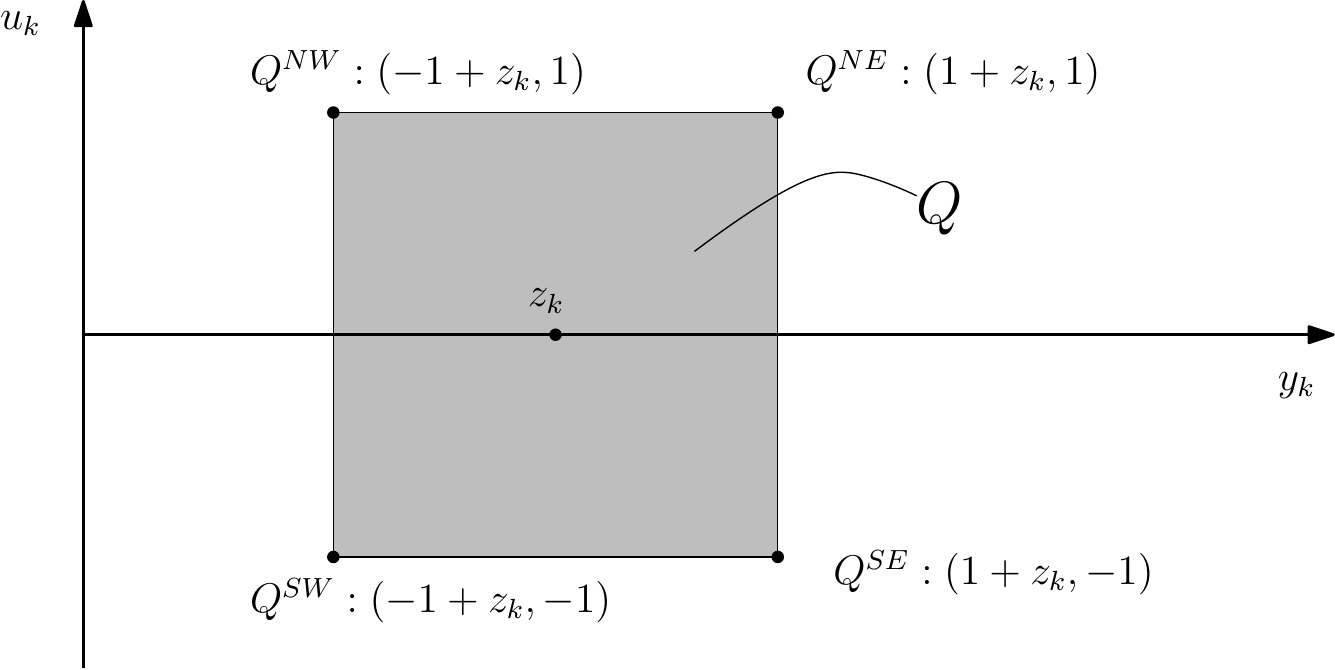} \caption{The square $Q$.}
	\label{Q}%
\end{figure}

See Fig. \ref{Q}. If the plant is in the state $\mathbf{x}$ at time $k$ then $Q(z_k) \cap L\left(\mathbf{x}\right)$ contains the plant's allowable outputs and inputs at time $k$.

We will also find it useful to define dual lines. Associated with any $\mathbf{f}\in \mathbb{R}^{m}$ define in the
$y^{\ast}_{k}u^{\ast}_{k}$-plane the \textit{dual line} $L^{\ast}\left(\mathbf{f}\right)$:

\begin{definition}
	$L^{\ast}\left(\mathbf{f}\right)=\{(y^*_k,u^*_k):n_{m+1}y^{\ast}_{k}+d_{m+1}u^{\ast}_{k}=-\left(\mathbf{f}\right)_1\}.$
\end{definition}
Suppose the dual plant is in the state $\mathbf{f}$ at time $k$. Then 
the scalars $y^*_k$ and $u^*_k$, the input and output of the dual plant, are constrained to lie on the line $L^{\ast}\left(  \mathbf{f}\right)$ by (\ref{adjointstate}) and (\ref{ABCDstardef2}).

\section{The primal estimation program and its dual} \label{sectprimaldual}

In this section we set up primal and dual optimisation programs. Their decision variables are the input and output signals of the primal and dual systems.
\subsection{Primal}
At any time $k > 2m$, and for any $\mathbf{f}^+\in\mathbb{R}^{m},$ consider the program 
\begin{equation} 
\max\limits_{\mathbf{x}^+%
	\in S_{k+1}}\left\langle \mathbf{f}^+,\mathbf{x}^+\right\rangle. 
\label{supequiv}
\end{equation}
It has optimal value $h_{S_{k+1}}(\mathbf{f}^+),$ the support function of $S_{k+1}$ evaluated at $\mathbf{f}^+$.

Writing out the constraints explicitly in terms of the output and input signals up to time $k$, namely $(\mathbf{y},\mathbf{u})=\left(  y_{1:k},u_{1:k}\right)$, by (\ref{conu_indconf}) the  program (\ref{supequiv}) can be equivalently expressed as: 
\[%
\begin{array}
[c]{c}%
\mathcal{P}_{z_{1:k}}(\mathbf{f}^+):\text{ \ \ }\max
\limits_{\mathbf{y},\mathbf{u}}\left\langle \mathbf{f}^+,\mathbf{x}_{k+1}\left(\mathbf{y},\mathbf{u}\right)%
\right\rangle \\
\text{subject to }\\
\mathbf{D}y_{1:k}-\mathbf{N}u_{1:k}=\left[
\begin{array}
[c]{c}%
\mathbf{B}_\mathrm{T}\mathbf{x}_{0} \\
\mathbf{0} 
\end{array} 
\right],  \\

\left\vert u_j\right\vert \leq1$ and $\left\vert y_j-z_j\right\vert \leq1$ for $j=1,\ldots, k.
\end{array} 
\]

From now on we will always assume that the constraints for the primal program are consistent, which is equivalent to saying that all of the uncertainty sets up to time $k$ are non-empty.

The following proposition is a consequence of the definitions and straightforward convexity arguments.

\begin{proposition}
	\label{propconeargmax} For any $\mathbf{x}\in S_{k}$ there exists $(\mathbf{y},\mathbf{u})$ feasible for $\mathcal{P}_{z_{1:k-1}}(\mathbf{\cdot})$ for which $\mathbf{x}=\mathbf{x}_k\left( \mathbf{y},\mathbf{u}\right)$. For any such  $(\mathbf{y},\mathbf{u})$, and any $\mathbf{f}\in\mathbb{R}^{m}$, there holds 
	
	\begin{align*}
	(\mathbf{y},\mathbf{u}) \in\arg\max\mathcal{P}_{z_{1:k-1}}(\mathbf{f})  &  \Leftrightarrow
	h_{S_{k}}(\mathbf{f})=\left\langle \mathbf{f},\mathbf{x} \right\rangle \\
	&  \Leftrightarrow
	\mathbf{f}\in \mathcal{N}_{S_{k}}\left(\mathbf{x} \right). \\
	\end{align*}
\end{proposition}

\subsection{Dual}
Let $\mathbf{f}^+\in\mathbb{R}^{m}$ be given. It will be shown in the next section that the dual of $\mathcal{P}_{z_{1:k}}(\mathbf{f}^+)$ is the program defined as follows.
\[%
\begin{array}
[c]{c}%
\mathcal{D}_{z_{1:k}}(\mathbf{f}^+):\text{ \ \ }\min\limits_{\mathbf{y}%
	^{\ast},\mathbf{u}^{\ast}}\left\{  \left\Vert \mathbf{y}^{\ast}\right\Vert
_{1}+\left\Vert \mathbf{u}^{\ast}\right\Vert _{1}+\left\langle \mathbf{y}^{\ast
},\mathbf{z}\right\rangle +\left\langle \mathbf{x}_{0}^{\ast}\left(\mathbf{y}
^{\ast},\mathbf{u}^{\ast}\right),\mathbf{x}
_{0}\right\rangle \right\} \\
\text{subject to }\\
\mathbf{N}^\mathrm{T}y^{\ast}_{1:k}+\mathbf{D}^\mathrm{T}u^{\ast}_{1:k}=\left[
\begin{array}
[c]{c}%
\mathbf{0}  \\
\mathbf{f}^+
\end{array}
\right].
\end{array}
\]
The decision variables, $(\mathbf{y}^{\ast},\mathbf{u}^{\ast}):=\left(  y_{1:k}^{\ast} ,u_{1:k}^{\ast}\right) $, are the inputs and outputs up to time $k$ of the dual system. The matrix $\mathbf{D}^\mathrm{T}$ $(\mathbf{N}^\mathrm{T})$ denotes the transpose of $\mathbf{D}$ $(\mathbf{N}),$ so the bottom right hand corner $m$ by $m$ submatrix of $\mathbf{D}^\mathrm{T}$ $(\mathbf{N}^\mathrm{T})$ is the transpose of $\mathbf{D}_\mathrm{L}\left(\mathbf{N}_\mathrm{L}\right).$ Thus, by (\ref{state_def_star}), the last $m$ of the constraint equations state that the decision variables are constrained by $\mathbf{x}^{\ast}_{k}\left(\mathbf{y}^{\ast},\mathbf{u}^{\ast}\right)=\mathbf{f}^+.$
\subsection{Alignment}
We define a relation between the inputs and outputs of the primal and dual systems. At optimality the primal and dual signals will be related through, in linear programming terminology, complementary slackness. The particular form this relationship takes in our setup is termed alignment, and is defined next. 
\begin{definition}
\label{def_align}The scalar pair $\left(  y_k,u_k\right)$ is said to be \textit{aligned}
	with $\left(  y_k^{\ast},u_k^{\ast}\right)  $ if
	\begin{equation*}%
	\begin{array}
	[c]{c}%
	u_k^{\ast}>0\implies u_k=1\\
	u_k^{\ast}<0\implies u_k=-1\\
	\left\vert u_k\right\vert <1\implies u_k^{\ast}=0
	\end{array}
	\end{equation*}
	and%
	\begin{equation*}%
	\begin{array}
	[c]{c}%
	y_k^{\ast}>0\implies y_k=1+z_{k}\\
	y_k^{\ast}<0\implies y_k=-1+z_{k}\\
	\left\vert y_k-z_{k}\right\vert <1\implies y_k^{\ast}=0.
	\end{array}
	\end{equation*}
\end{definition}
This definition can be extended in a natural way to alignment between pairs of
sequences. Thus the vector pair $\left( \mathbf{y},\mathbf{u}\right)  $ is
aligned with the pair $\left(  \mathbf{y}^{\ast},\mathbf{u}^{\ast}\right)  $ if, for
all $k$, $\left(  y_{k},u_{k}\right)  $ is aligned with $\left(
y_{k}^{\ast},u_{k}^{\ast}\right)  .$
Alignment between points in the $y_ku_k$ and $y_k^{\ast}u_k^{\ast}$-planes can be readily visualised as follows.

\begin{figure}[ptb]
	\centering
	\includegraphics[scale=0.7]{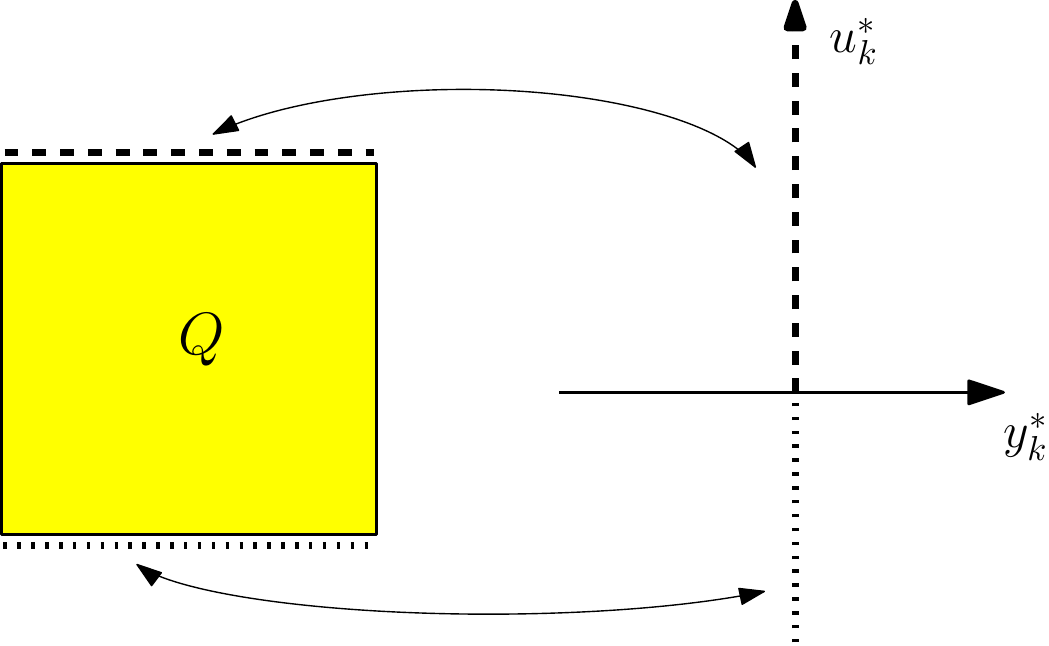} \caption{Alignment between the top (bottom) side of Q and the positive (negative) $u_k^{\ast}$ axis.}
	\label{align_TB}%
\end{figure}

\begin{figure}[ptb]
	\centering
	\includegraphics[scale=0.55]{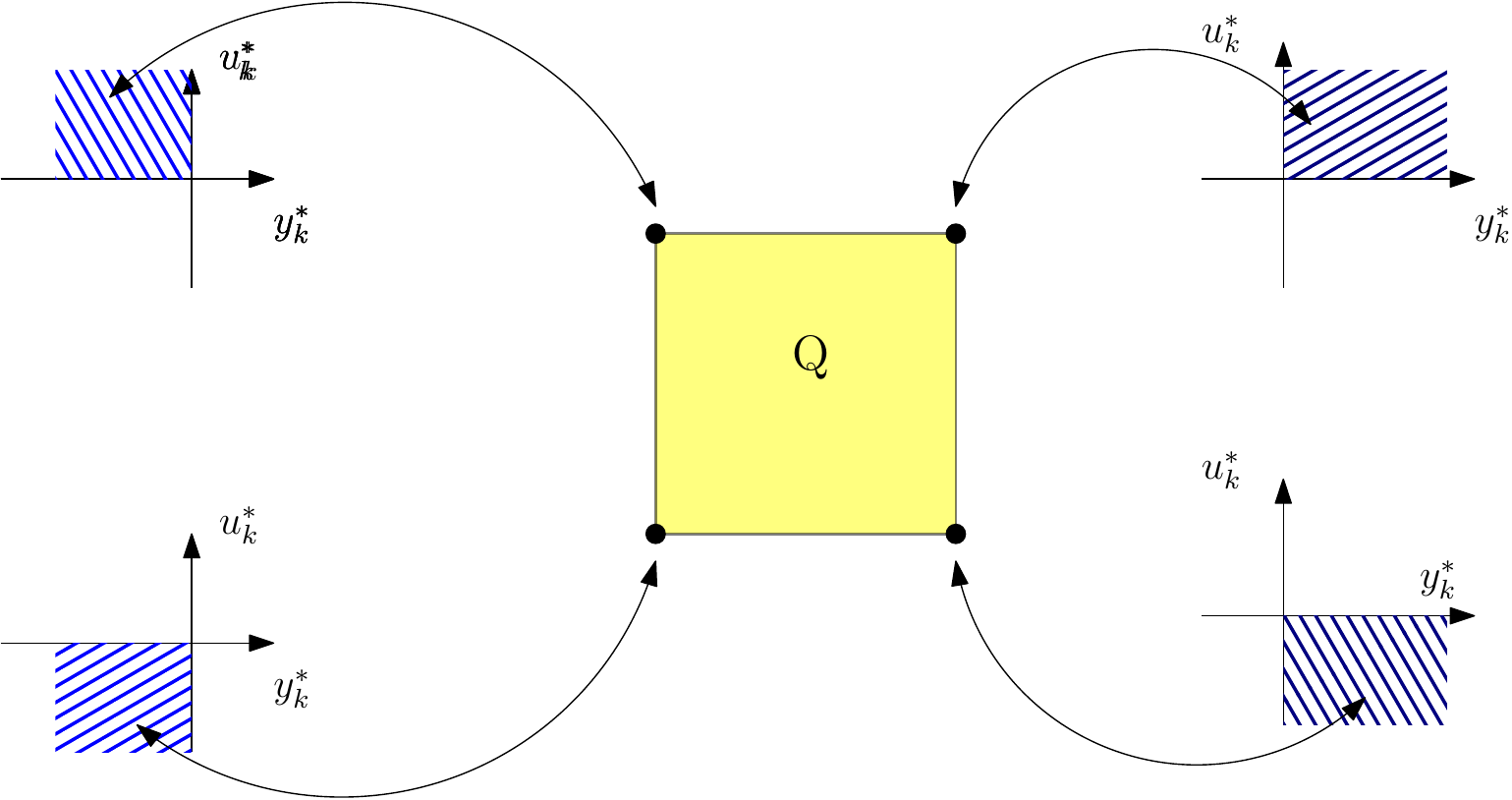} \caption{Alignment between the corners of Q and the four quadrants in the $y_k^{\ast}u_k^{\ast}$-plane.}
	\label{align_corners}%
\end{figure}

By Definition (\ref{def_align}) each point $(y_k,u_k)$ belonging to the top (bottom) side of $Q$ is aligned with every point on the positive (negative) $u_k^{\ast}$ axis in the $y_k^{\ast}u_k^{\ast}$-plane. See Fig. \ref{align_TB}. Similarly, each point $(y_k,u_k)$ belonging to the right (left) side of $Q$ is aligned with every point on the positive (negative) $y_k^{\ast}$ axis in the $y_k^{\ast}u_k^{\ast}$-plane. 
The corner $Q^{\mathrm{NE}} (Q^{\mathrm{NW}},Q^{\mathrm{SW}},Q^{\mathrm{SE}})$ of $Q$ is aligned with all points in the first (respectively, second, third, fourth) quadrant of the $y_k^{\ast}u_k^{\ast}$-plane. See Fig. \ref{align_corners}.

Finally, all points in Q are aligned with the origin in the $y_k^{\ast}u_k^{\ast}$-plane.
\subsection{Fundamental theorem}

A formal statement of the duality relating $\mathcal{P}_{z_{1:k}%
}(\mathbf{f}^+)$ to $\mathcal{D}_{z_{1:k}}(\mathbf{f}^+)$
is now presented. It is the basis for the proofs of our main theorems. Amongst other things, it justifies the interpretation of $\mathbf{x}^*_{k+1}\left(\mathbf{y}^*,\mathbf{u}^*\right)$ as the direction vector $\mathbf{f}^+$, the argument of the support function of $S_{k+1}$, see (\ref{supequiv}). It is worth pointing out that the structurally elegant form manifest in this result is not apparent in a routine application of duality to $\mathcal{P}_{z_{1:k}}(\mathbf{f}^+)$. Observe, for example, that for the program $\mathcal{D}_{z_{1:k}}(\mathbf{f}^+)$ the initial state
is free, and the terminal state is fixed, at $\mathbf{f}^+$. For the
program $\mathcal{P}_{z_{1:k}}(\mathbf{f}^+)$ the initial state
is fixed, at $\mathbf{x}_{0}$, and the terminal state is free. Some finesse is required in the construction of the dual variables and the dual cost function. Any dual of $\mathcal{P}_{z_{1:k}}(\mathbf{f}^+)$ will be equivalent to $\mathcal{D}_{z_{1:k}}(\mathbf{f}^+)$, but the fact that duality can be used to prove Theorem \ref{main1} only becomes apparent when the dual is expressed in the form $\mathcal{D}_{z_{1:k}}(\mathbf{f}^+)$.

\begin{theorem}
	\label{main}Suppose the set $S_{k+1}$ is non-empty. Then the optimal
	values of $\mathcal{P}_{z_{1:k}}(\mathbf{f}^+)$ and
	$\mathcal{D}_{z_{1:k}}(\mathbf{f}^+)$ are finite and equal. Furthermore,
	if $(\mathbf{y},\mathbf{u})$ and $(\mathbf{y}^{\ast},\mathbf{u}^{\ast})$ are
	feasible for $\mathcal{P}_{z_{1:k}}(\mathbf{f}^+)$ and
	$\mathcal{D}_{z_{1:k}}(\mathbf{f}^+),$ respectively, then a necessary and
	sufficient condition that they both be optimal is that they be aligned.
\end{theorem}

proof

The proof in outline follows standard linear programming arguments, although some non-routine manipulations involving the
	Gohberg-Semencul formula (\ref{def:BTmatrix}) are also required. Details are in
	the Appendix.

\begin{remark}
	It can be shown that, under the standing assumption that $\hat{\mathbf{n}}$	and $\hat{\mathbf{d}}$ are coprime, $\mathcal{D}_{z_{1:k}}(\mathbf{f}^+)$ always has a feasible solution, and has unbounded negative cost if $S_{k+1}$ is empty. 
\end{remark}

We now define a set $M$ of pairs of scalars, the first element of a pair being a possible input to the primal plant at time $k$, and the second a special input, related through alignment, to the dual plant.   

\begin{definition}
		\label{defM}Given $\mathbf{x}\in S_{k}$, $\mathbf{f}\in \mathcal{N}_{S_{k}}(\mathbf{x})$  and $z_{k},$ the set $M\left(
	\mathbf{x},\mathbf{f},z_{k}\right)$ is the set of scalar pairs $\left(  u_k,y_k^{\ast}\right)$ which satisfy
	\begin{enumerate}
		\item there exists $(y_k,u_k) \in L\left(  \mathbf{x}\right) \cap Q,$ and 
		\item  $\left(y_k,u_k\right)$  is aligned with $\left(
		y_k^{\ast},u_k^{\ast}\right)$,  where
		$(y_k^{\ast},u_k^{\ast}) \in L^{\ast}\left(  \mathbf{f}\right).$
	\end{enumerate}
\end{definition}

\begin{figure}[ptb]
	\begin{center} 
		\includegraphics[scale=0.5]{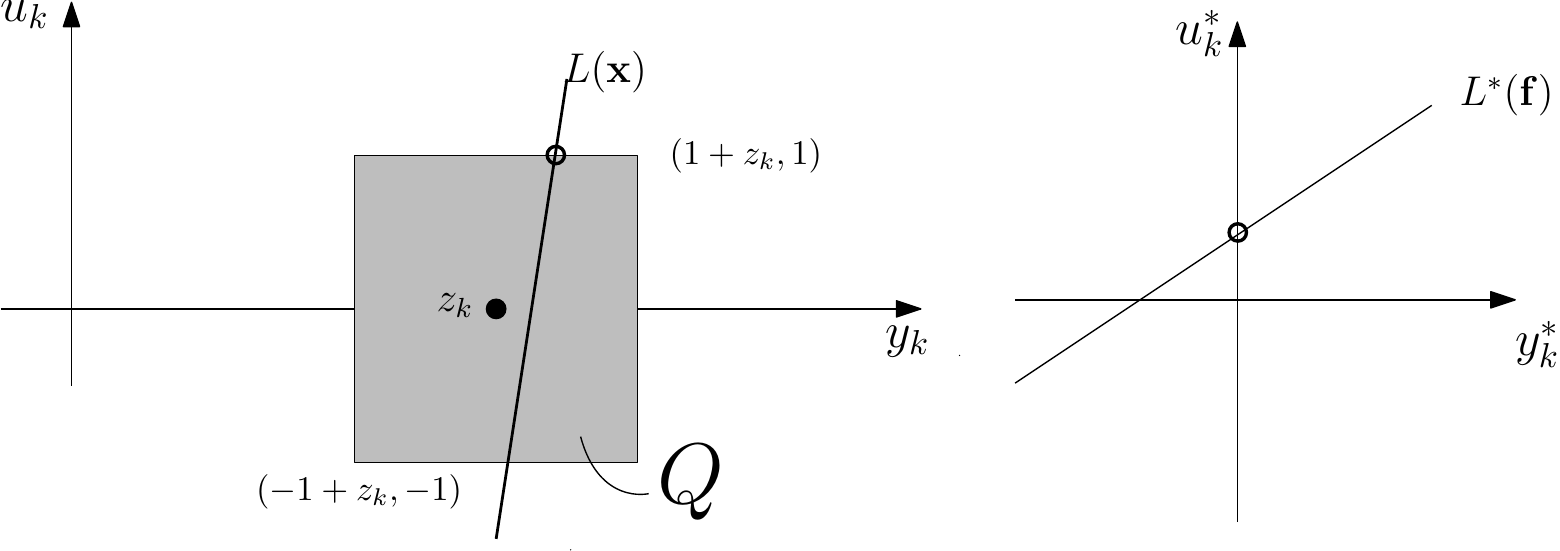} \caption{The circled points on the lines $L(\mathbf{x})$ and $L^{\ast}(\mathbf{f})$ are aligned, so $M\left(
			\mathbf{x},\mathbf{f},z_{k}\right)=\{(1,0)\}$ because at the aligned points $u_k=1$ and $y^*_k=0$.} \label{Q_L_intersect}%
	\end{center}

\end{figure}

Finding $M$ is computationally very simple, requiring merely the intersection of lines in the plane, and checking alignment. For example, in Fig. \ref{Q_L_intersect} alignment for points in $L\left(  \mathbf{x}\right) \cap Q$ occurs solely between the two circled points, so $M$ contains the singleton pair $(u_k,y^*_k)=(1,0)$.

\section{Main results}

\begin{figure*}[ptb]
	\centering
	\includegraphics[scale=0.9]{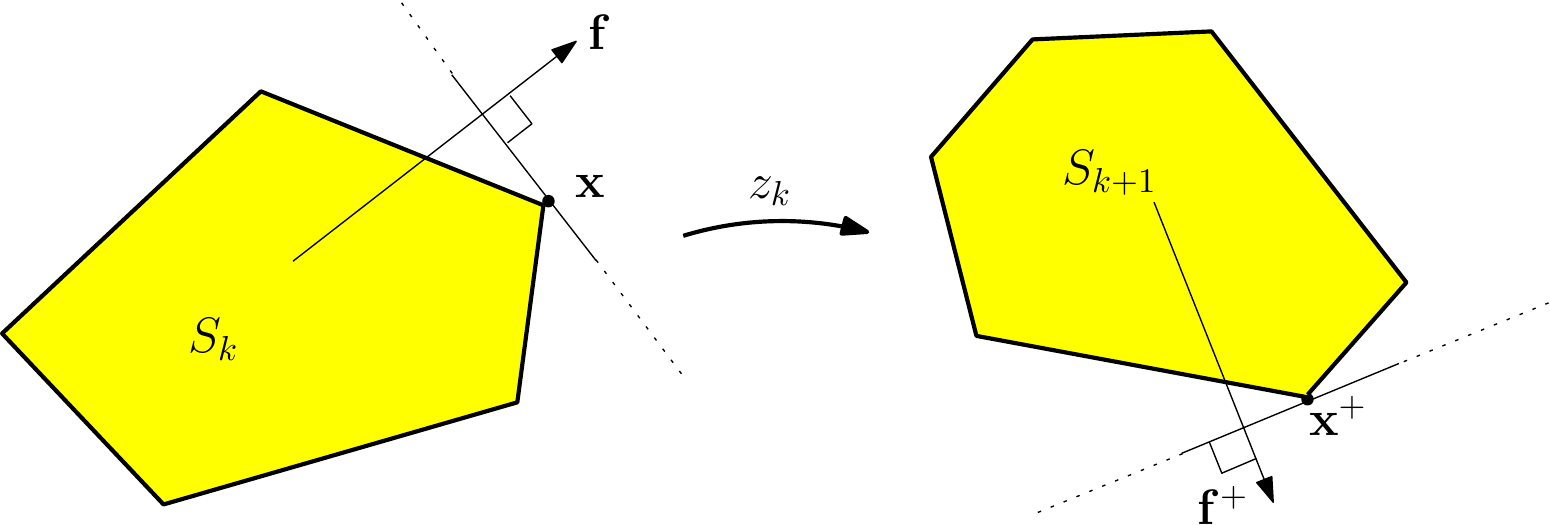} \caption{The vectors in Theorem \ref{main1}. The state
		$\mathbf{x}\in\partial S_{k}$ and direction $\mathbf{f}\in
		\mathcal{N}_{S_{k}}(\mathbf{x})$ are propagated to $\mathbf{x^+}%
		\in\partial S_{k+1}$ and $\mathbf{f^+}\in \mathcal{N}_{S_{k+1}}(\mathbf{x^+}%
		)$ by the measurement $z_{k}.$}%
	\label{basic_conv}%
\end{figure*}

We now present Theorems \ref{main1} and \ref{main2}. See Fig. \ref{basic_conv} for a geometric depiction of the vectors in these theorems for the special case where $\mathbf{x}$ and $\mathbf{x}^+$ are both boundary points. Proofs are in the Appendix.
\vspace{5pt}

\begin{theorem}
	\label{main1}Suppose $\mathbf{x}\in S_{k}$ and
	$\mathbf{f}\in \mathcal{N}_{S_{k}}(\mathbf{x})$. Then $\mathbf{x^+}=\mathbf{A}\mathbf{x}+\mathbf{B}u_{k} \in S_{k+1}$ and 
	$\mathbf{f^+}
	=\mathbf{A}^{\ast}\mathbf{f}+\mathbf{B}^{\ast}y_{k}^{\ast
	}\in \mathcal{N}_{S_{k+1}}(\mathbf{x^+})$ if and only if 	$\left(
	u_{k},y_{k}^{\ast}\right)  \in M\left( \mathbf{x}, \mathbf{f},z_{k}\right)$.
\end{theorem}

For given $\mathbf{x} \in S_{k}$ and $\mathbf{f}\in \mathcal{N}_{S_{k}}(\mathbf{x})$ Theorem \ref{main1} furnishes at least one successor $\mathbf{x}^+$ to $\mathbf{x}$, and one vector $\mathbf{f}^+ \in \mathcal{N}_{S_{k+1}}\left( \mathbf{x}^+\right)$, if $ M\left( \mathbf{x}, \mathbf{f},z_{k}\right)$ is non-empty.  It is not yet clear, however, that \textit{every} point of $S_{k+1}$ can be found through the process of applying Theorem \ref{main1} to some $\mathbf{x} \in S_{k}$ and some $\mathbf{f}\in \mathcal{N}_{S_{k}}(\mathbf{x})$.

The companion result Theorem \ref{main2}, given below, shows that any boundary point $\mathbf{x^+} \in S_{k+1}$, and any direction in the normal cone of $\mathbf{x^+}$, are attainable from \textit{any} precursor $\mathbf{x}$ of $\mathbf{x^+}$  and \textit{some} direction in the normal cone of $\mathbf{x}$. 

\begin{theorem}
	\label{main2}	
	Select any $\mathbf{x^+}\in S_{k+1}$, any $\mathbf{f^+} \in \mathcal{N}_{S_{k+1}}\left( \mathbf{x^+}\right)$ and any precursor $\mathbf{x}$ of $\mathbf{x^+}$. There exists $\mathbf{f} \in \mathcal{N}_{S_{k}}\left( \mathbf{x}\right)$ and $\left(
	u_{k},y_{k}^{\ast}\right)  \in M\left( \mathbf{x},  \mathbf{f},z_{k}\right)$
	for which $\mathbf{x^+} =\mathbf{A}\mathbf{x}+\mathbf{B}u_{k}$ and $\mathbf{f^+}=\mathbf{A}^{\ast}\mathbf{f}+\mathbf{B}^{\ast}y_{k}^{\ast}$.
\end{theorem}

The proof relies on a dynamic programming style argument and is given in the Appendix.


Theorem \ref{main1} shows how any point in the uncertainty set together with any direction in its normal cone can be propagated efficiently and recursively. Theorem \ref{main2} shows that every point in $S_{k+1}$, and every direction in its normal cone, is the result of applying Theorem \ref{main1} to some point in $S_k$, along with a direction in its normal cone. It follows that the vertices and facet directions of $S_{k+1}$ can be determined by propagating the vertices and facet directions of $S_k$. Implementation of this simple idea is somewhat involved and space limitations preclude giving more than a sketch of the algorithm that achieves this. However, code implementing our algorithm for plants whose primal and dual lines have positive slopes, called \texttt{uncertaintyset.m}, is available on the link below \footnote{Go to \textrm{https://mathworks.com/matlabcentral/fileexchange}}. Also available on this link is the code \texttt{tcomp.m}. It generates plants and measurements randomly, and compares the performance of three algorithms. 

For the special case of plants with a lag, that is $n_1=0$, the algorithmic details are considerably simplified. Some simulations for this case are in \cite{HILL-CDC2017}.

We sketch how the updating is done for general linear time-invariant plants. Suppose $F$ is a facet of $S_k$, and the vertices and direction of $F$ are, respectively, $\mathbf{v}_j$ and $\mathbf{f}$. A simple case, which nonetheless illustrates how Theorem \ref{main1} is useful in the propagation of facets, occurs when all primal lines $L\left(\mathbf{v}_j\right)$ pass through the interior of the same side of $Q$, and the dual line $L^{\ast}\left(\mathbf{f}\right)$ does not pass through the origin. Suppose, for example, that all $L\left(\mathbf{v}_j\right)$ pass through the interior of the top side of $Q$, and $L^{\ast}\left(\mathbf{f}\right)$ intersects the positive $u_k^{\ast}$ axis in the $y_k^{\ast}u_k^{\ast}$-plane, see Fig. \ref{align_TB}. Then it follows from Theorem \ref{main1} that there is a facet of $S_{k+1}$ with vertices $\mathbf{A}\mathbf{v}_j+\mathbf{B}$, and direction $\mathbf{A}^{\ast}\mathbf{f}.$  Many facets of $S_{k+1}$ can be identified very quickly in similar fashion.

There are two cases in which complications to the simple propagation described in the previous paragraph can arise. The first is when primal lines of states in $S_k$ intersect a corner point of $Q$, and the second is when dual lines pass through the origin in the $y_k^{\ast}u_k^{\ast}$-plane. In the first case vertices of $S_{k+1}$ arise which are not affine images of vertices of $S_k$; instead they are affine images of vertices formed from the intersection of hyperplanes with $S_k$. In the second case new possibilities for facet propagation are opened up because
all points in $Q$ are aligned with the origin in the $y_k^{\ast}u_k^{\ast}$-plane; in fact ridges in $S_k$ can propagate to facets in $S_{k+1}$. There is insufficient space to give details here, but we claim that our algorithm updates both vertices and facets. The only computationally intensive tasks are intersecting at most four hyperplanes with $S_k$, and the intersecting all pairs of facets of $S_k$ whose directions have first components of opposite sign. It is a consequence of Theorems \ref{main1} and \ref{main2} that the only intersections needed are those which are guaranteed to produce propagated vertices and/or facets. No time is wasted in calculating, and then discarding, redundant inequality constraints, as occurs in current exact methods.

\section{Numerical simulations}

The accompanying code is written for the special case of plants whose primal and dual lines both have a positive slope, that is $n_1 > 0$ and $n_{m+1}d_{m+1} < 0$. Other cases are not any more difficult; they just require modifications to the coding. We present results obtained by running this code for randomly selected stable plants. We give some comparisons of computation time for our algorithm, denoted FV (facets-vertices), Fourier-Motzkin and plp, which is an acronym for parametric linear programming. These last two are commonly used in applications requiring recursive determination of uncertainty sets. We make use of the multi-parametric toolbox \cite{mpt} for their implementation.

\begin{figure}[ptb]
	\centering
	\includegraphics[scale=0.6]{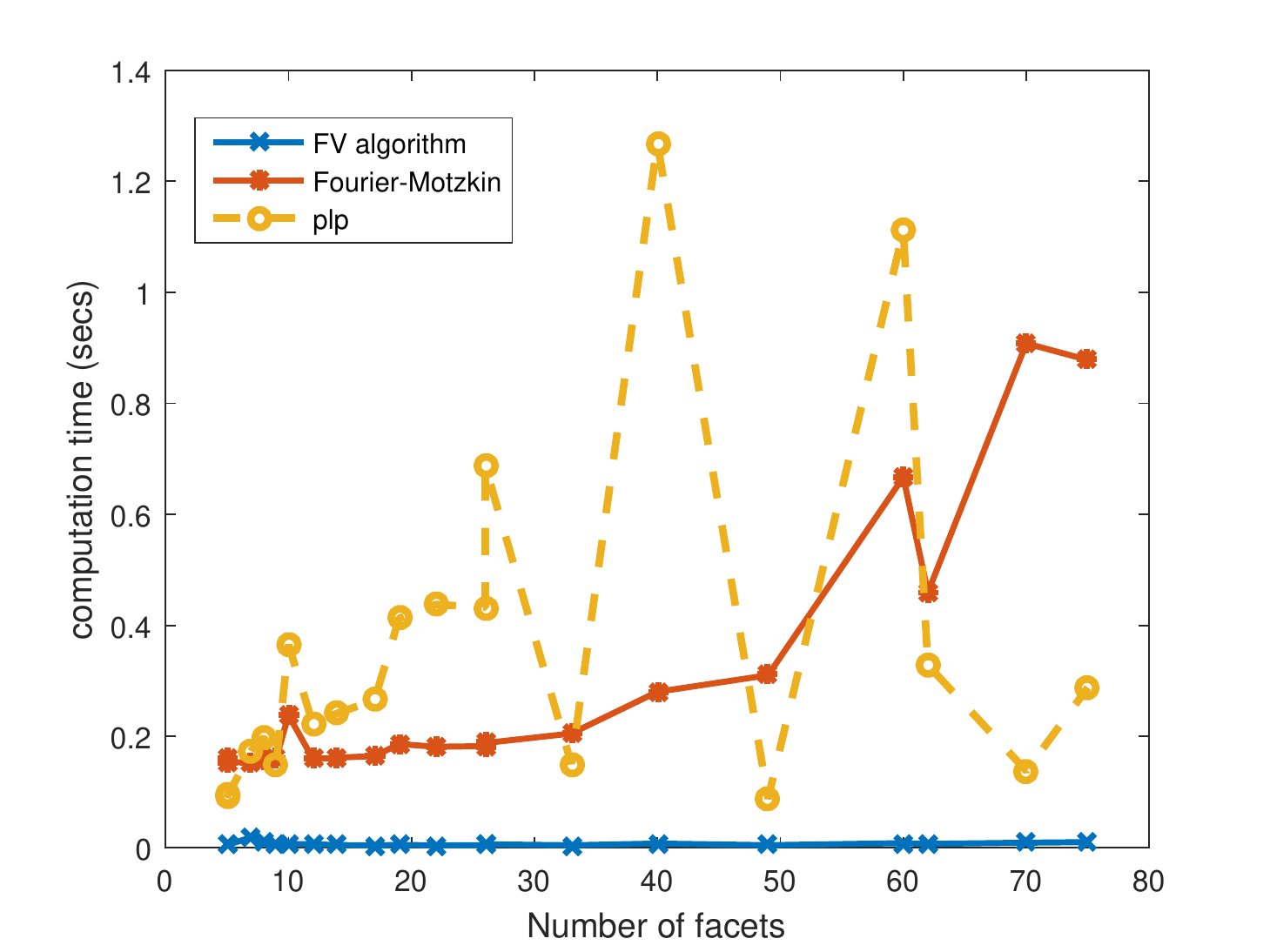} \caption{The time taken to update $S_k$ versus the number of facets of $S_k$ for a third order plant }%
	\label{timecomp_3rdorder}%
\end{figure}

\begin{figure}[ptb]
	\centering
	\includegraphics[scale=0.6]{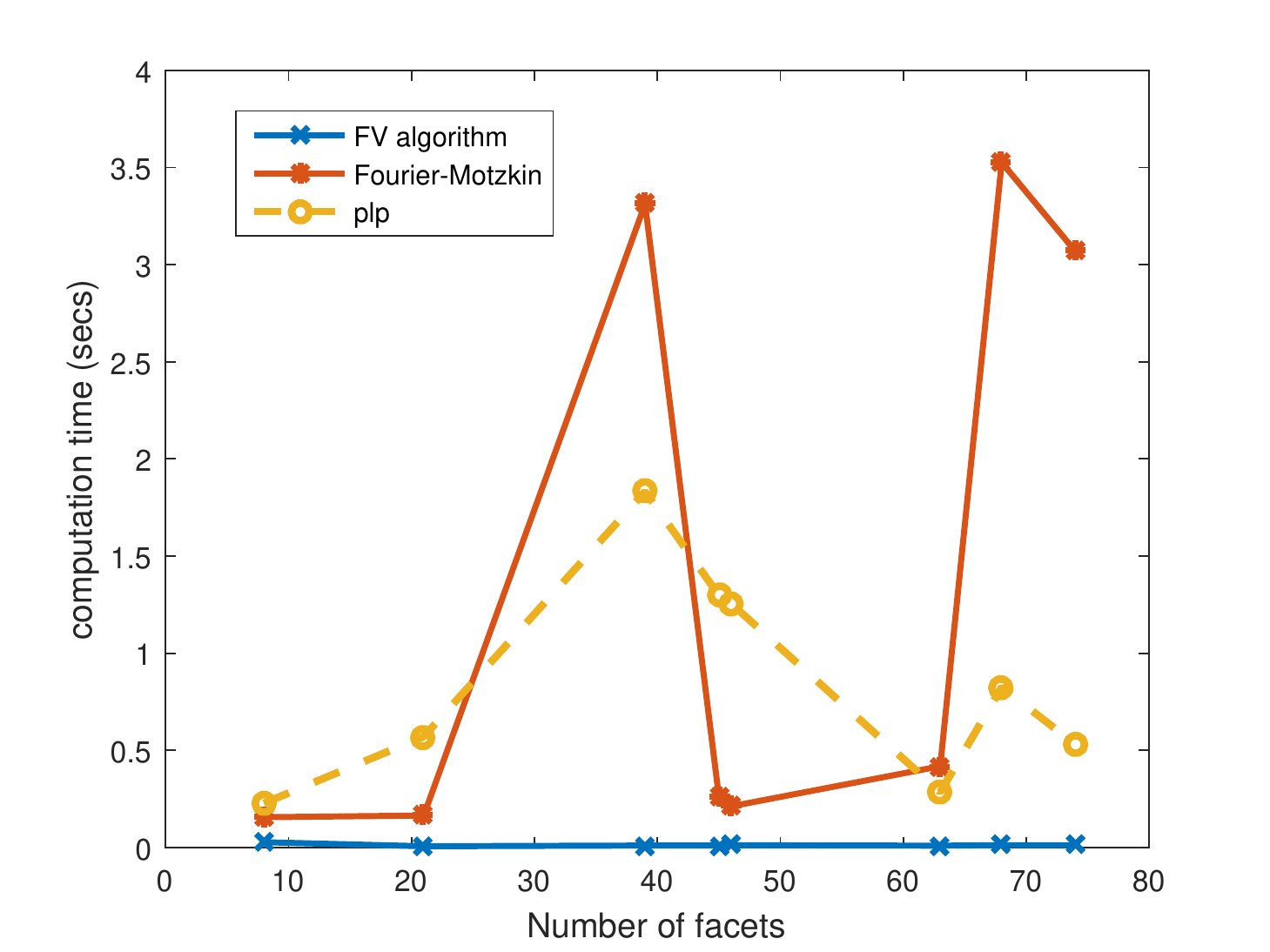} \caption{The time taken to update $S_k$ versus the number of facets of $S_k$ for a fourth order plant }%
	\label{comptime_4thorder}%
\end{figure}

Fig. \ref{timecomp_3rdorder} shows results for a third order randomly selected plant and random measurements, with the initial uncertainty set being a simplex with four vertices and four facets. Along the horizontal axis the number of facets of the polytopic uncertainty set is displayed. Along the vertical axis is the computation time required for one and the same update by the three algorithms. The FV algorithm is seen to be the fastest. A similar pattern is seen in Fig. \ref{comptime_4thorder}, which simulates uncertainty set propagation for a randomly selected fourth order plant, again with random measurements. For this example the average ratio of update times for plp and FV was $80.5$; for Fourier-Motzkin and FV it was $123$. These ratios increase with increasing complexity of the uncertainty set.

\begin{figure}[ptb]
	\centering
	\includegraphics[scale=0.6]{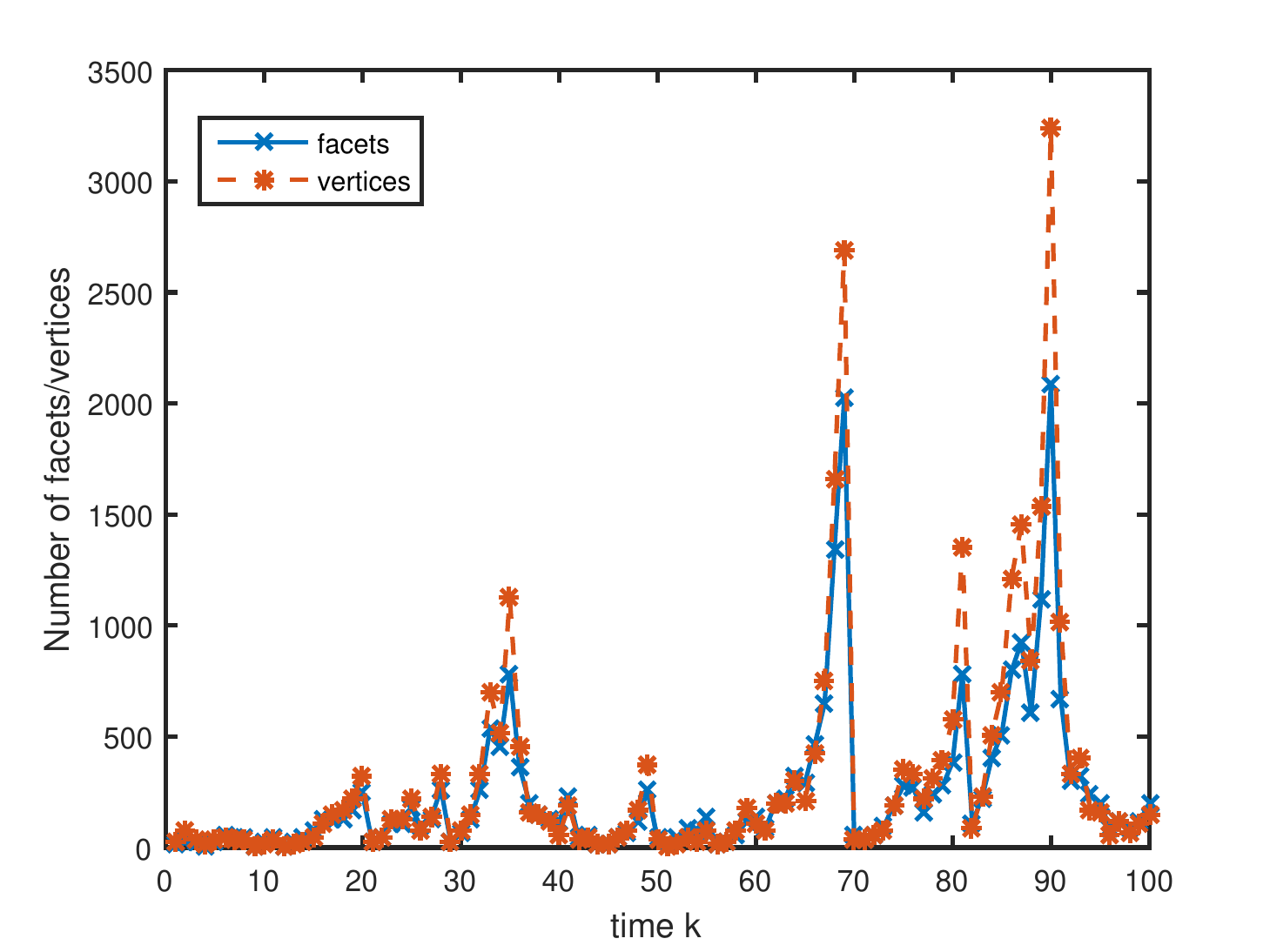} \caption{The number of facets and vertices of $S_k$ for a fifth order plant and random measurements. }%
	\label{vertices_growth}%
\end{figure}

The update computation time for Fourier-Motzkin and plp becomes prohibitive when the number of facets exceeds a couple of hundred. The computation time also increases when using the FV algorithm, but we can continue updating for much longer. The growth in the number of vertices and facets for a randomly selected fifth order plant and one hundred time instants using the FV algorithm is depicted in Fig. \ref{vertices_growth}. The  complexity of the uncertainty set varies quite dramatically with time, but does not appear to be consistently increasing. For plants of order higher than five there is typically a very rapid increase in the number of facets and vertices. 

\section{CONCLUSIONS}

We have introduced two theorems which describe completely the evolution of state uncertainty sets. When implemented in code there appears to be a significant performance improvement over existing exact methods.

More work is needed to extend the results to time-varying linear plants, and to multivariable systems. There does not seem to be any reason to believe this will not be possible. For time-varying plants the primal lines will no longer be parallel as time evolves, and the same statement holds true for dual lines, but checking alignment is not any more difficult. For multivariable plants, alignment can be defined between vector inputs and outputs, although the computations will be more involved. The use of state and associated normal cone propagation in multivariable systems is a topic for future research.

\bibliographystyle{plain} 
\bibliography{acompat,Ref}

\appendix
\section{Proof of Theorem \ref{main}}
After expressing the program $\mathcal{D}_{z_{1:k}}(\mathbf{f}^+)$ as an
equivalent linear program, the standard duality result in asymmetric form
(\cite{LUEN-84} p. 86, 96) is used:%
\begin{equation}%
	\begin{array}
		[c]{cc}%
		\text{Dual} & \text{Primal}\\%
		\begin{array}
			[c]{c}%
			\min\mathbf{c}^\mathrm{T}\mathbf{x}\\
			\text{s. t. }\mathbf{A}\mathbf{x}=\mathbf{b}\\
			\mathbf{x}\geq0
		\end{array}
		&
		\begin{array}
			[c]{c}%
			\max\bm{\theta}^\mathrm{T}\mathbf{b}\\
			\text{s. t. }\mathbf{A}^\mathrm{T}\bm{\theta}\leq\mathbf{c}%
		\end{array}
		\text{,}%
	\end{array}
	\label{asymdual}%
\end{equation}
where \textit{complementary slackness} holds: Let $\mathbf{x}$ and
$\bm{\theta}$ be feasible solutions for the primal and dual problems,
respectively. A necessary and sufficient condition that they both be optimal
solutions is that for all $i$

i) $x_{i}>0\implies a_{i}^\mathrm{T}\bm{\theta}=c_{i}$ (where $a_{i}^\mathrm{T}$ is the $i^{\rm th}$
row of $\mathbf{A}^\mathrm{T}$)

ii) $x_{i}=0\Leftarrow a_{i}^\mathrm{T}\bm{\theta}<c_{i}.$

Note that we have swapped the labels for primal and dual from that given in \cite{LUEN-84}, because it is the estimation program that we regard as the primal problem, and the estimation program is naturally expressed in the form of the maximization in (\ref{asymdual}).  Also the use of the symbol $\mathbf{x}$ for the dual decision variable
in (\ref{asymdual}) is different from the use of the symbols $\mathbf{x}_{0}$,
$\mathbf{x}_{0}^{\ast}$ and $\mathbf{x}_{k}$, which
retain their meanings given in the body of the paper.

The program $\mathcal{D}_{z_{1:k}}(\mathbf{f})$ has a convex
piecewise linear cost function and linear constraints. There is a standard
procedure, which we now follow, for converting such a program to an equivalent
linear programming problem. Introduce new non-negative $k$-dimensional column
vectors $\mathbf{u}^{\ast+},\mathbf{u}^{\ast-},\mathbf{y}^{\ast+}$ and
$\mathbf{y}^{\ast-}$, and put $u_{j}^{\ast}=u_{j}^{\ast+}-u_{j}^{\ast-}$ and
$y_{j}^{\ast}=y_{j}^{\ast+}-y_{j}^{\ast-}.$ At optimality at least one of
$u_{j}^{\ast+},u_{j}^{\ast-}$, and at least one of $y_{j}^{\ast+},y_{j}%
^{\ast-}$, will be zero, so $\left\vert u_{j}^{\ast}\right\vert =u_{j}^{\ast
	+}+u_{j}^{\ast-}$ and $\left\vert y_{j}^{\ast}\right\vert =y_{j}^{\ast+}%
+y_{j}^{\ast-}$. Since $\left\langle \mathbf{x}_{0}^{\ast},\mathbf{x}%
_{0}\right\rangle =-\mathbf{x}_{0}^\mathrm{T}\left[  \mathbf{N}_\mathrm{U}^\mathrm{T}y_{1:m}^{\ast}%
+\mathbf{D}_\mathrm{U}^\mathrm{T}u_{1:m}^{\ast}\right]  $, the cost function for
$\mathcal{D}_{z_{1:k}}(\mathbf{f}^+)$, namely $\left\Vert y^{\ast
}\right\Vert _{1}+\left\Vert u^{\ast}\right\Vert _{1}+\left\langle
y_{1:k}^{\ast},z_{1:k}\right\rangle +\left\langle \mathbf{x}_{0}^{\ast
},\mathbf{x}_{0}\right\rangle =:J_{\mathrm{dual}},$ can be written as
\[
J_{\mathrm{dual}}=\left[  \boldsymbol{1}_{4k}+\bm{\delta}+\bm{\gamma}\right]
\left[
\begin{array}
[c]{c}%
\mathbf{y}^{\ast+}\\
\mathbf{y}^{\ast-}\\
\mathbf{u}^{\ast+}\\
\mathbf{u}^{\ast-}%
\end{array}
\right]
\]
where$\ \boldsymbol{1}_{4k}$ denotes a $4k-$dimensional row vector of ones,
the row vector $\bm{\delta}$ is defined to be
\begin{align*}
\left[
\begin{array}
{cccccccc}%
-\mathbf{x}_{0}^\mathrm{T}\mathbf{N}_\mathrm{U}^\mathrm{T} & \mathbf{0}_{k-m} & \mathbf{x}_{0}^\mathrm{T}\mathbf{N}_\mathrm{U}^\mathrm{T}
& \mathbf{0}_{k-m} & -\mathbf{x}_{0}^\mathrm{T}\mathbf{D}_\mathrm{U}^\mathrm{T} & \mathbf{0}_{k-m} &
\mathbf{x}_{0}^\mathrm{T}\mathbf{D}_\mathrm{U}^\mathrm{T} & \mathbf{0}_{k-m}%
\end{array}
\right]
\end{align*}
where $\mathbf{0}_{k-m}$ denotes a $\left(  k-m\right)  $-dimensional row
vector of zeros,
and %
\begin{align*}
	\bm{\gamma}  :=\left[
	\begin{array}
		{ccc}%
		z_{1:k}^\mathrm{T} & -z_{1:k}^\mathrm{T} & \mathbf{0}_{2k}%
	\end{array}
	\right].
\end{align*}
The constraints for the program $\mathcal{D}_{z_{1:k}}(\mathbf{f}^+)$\ in
terms of the new variables are%
\begin{align*}
	\left[
	\begin{array}
		[c]{cccc}%
		\mathbf{N}^\mathrm{T} & -\mathbf{N}^\mathrm{T} & \mathbf{D}^\mathrm{T} & -\mathbf{D}^\mathrm{T}%
	\end{array}
	\right]  \left[
	\begin{array}
		[c]{c}%
		\mathbf{y}^{\ast+}\\
		\mathbf{y}^{\ast-}\\
		\mathbf{u}^{\ast+}\\
		\mathbf{u}^{\ast-}%
	\end{array}
	\right]   &  =\left[
	\begin{array}
		[c]{c}%
		0\\
		\vdots\\
		0\\
		\mathbf{f}^+%
	\end{array}
	\right] \\
	y_{j}^{\ast+},y_{j}^{\ast-},u_{j}^{\ast+},u_{j}^{\ast-}  &  \geq0.
\end{align*}

In (\ref{asymdual}) put
\begin{align}
	A  &  =\left[
	\begin{array}
		[c]{cccc}%
		\mathbf{N}^\mathrm{T} & -\mathbf{N}^\mathrm{T} & \mathbf{D}^\mathrm{T} & -\mathbf{D}^\mathrm{T}%
	\end{array}
	\right] \nonumber 
	\\
	\mathbf{x}  &  =\left[
	\begin{array}
		[c]{cccc}%
		\mathbf{y}^{\ast+\mathrm{T}} & \mathbf{y}^{\ast-\mathrm{T}} & \mathbf{u}^{\ast+\mathrm{T}} &
		\mathbf{u}^{\ast-\mathrm{T}}%
	\end{array}
	\right]  ^\mathrm{T},\text{ }\bm{c}^\mathrm{T}=\boldsymbol{1}_{4k}+\bm{\delta}+\bm{\gamma
	}\label{notation_dual}\\
	\mathbf{b}  &  =\left[
	\begin{array}
		[c]{cccc}%
		0 & \ldots & 0 & \mathbf{f}^T%
	\end{array}
	\right]  ^\mathrm{T},\nonumber
\end{align}
so the program labelled Dual in the left column of (\ref{asymdual}) is equivalent to $\mathcal{D}_{z_{1:k}}(\mathbf{f}^+)$.

Then by (\ref{asymdual}) the dual of $\mathcal{D}%
_{z_{1:k}}(\mathbf{f}^+)$ is%
\[
\begin{array}
[c]{c}%
\bar{\mathcal{P}}(\mathbf{f}^+):\text{ \ \ }\max\limits_{\bm{\theta}\in{\mathbb{R}}^{k}}\left\langle \bm{\theta
}_{k-m+1:k},\mathbf{f}^+\right\rangle \\
\text{subject to }\\
\left[
\begin{array}
[c]{c}%
\mathbf{N}\\
-\mathbf{N}\\
\mathbf{D}\\
-\mathbf{D}%
\end{array}
\right]  \bm{\theta}\leq\left[  \boldsymbol{1}_{4}+\bm{\delta
}+\bm{\gamma}\right]  ^\mathrm{T}.
\end{array}
\]
The task from now on is to show that $\bar{\mathcal{P}}(\mathbf{f}^+)$ is equivalent to $\mathcal{P}_{z_{1:k}}(\mathbf{f}^+).$ The proof of this is done in two parts. In the first part we show that any  $\bm{\theta}$ feasible for $\bar{\mathcal{P}}(\mathbf{f}^+)$ corresponds to a feasible solution to $\mathcal{P}_{z_{1:k}}(\mathbf{f}^+)$ with the same cost. In the second part we show that any feasible solution to $\mathcal{P}_{z_{1:k}}(\mathbf{f}^+)$ corresponds to a feasible solution to $\bar{\mathcal{P}}(\mathbf{f}^+)$ with the same cost.

Before launching into details we mention that upper triangular Toeplitz matrices commute, so for example $\mathbf{N}_\mathrm{U}\mathbf{D}_\mathrm{U}=\mathbf{D}_\mathrm{U}\mathbf{N}_\mathrm{U}$. We shall frequently make use of this, as well as the similar comment which can be made for lower triangular Toeplitz matrices.

For the first part, take any $\bm{\theta}$ feasible for $\bar{\mathcal{P}}(\mathbf{f}^+)$, and put

\begin{equation}
	\mathbf{u}:=\mathbf{D}\bm{\theta}+\left[
	\begin{array}
		[c]{c}%
		\mathbf{D}_\mathrm{U}\mathbf{x}_{0}\\
		0
	\end{array}
	\right]  \text{and }\mathbf{y}:=\mathbf{N}\bm{\theta}+\left[
	\begin{array}
		[c]{c}%
		\mathbf{N}_\mathrm{U}\mathbf{x}_{0}\\
		0
	\end{array}
	\right]  \label{eq13}%
\end{equation}
which is equivalent to putting %
\begin{equation}
	\left[
	\begin{array}
		[c]{c}%
		\mathbf{y}\\
		-\mathbf{y}\\
		\mathbf{u}\\
		-\mathbf{u}%
	\end{array}
	\right]
	=
	\left[
	\begin{array}
		[c]{c}%
		\mathbf{N}\\
		-\mathbf{N}\\
		\mathbf{D}\\
		-\mathbf{D}%
	\end{array}
	\right]  \bm{\theta}-\bm{\delta}^\mathrm{T}. \label{eq14}%
\end{equation}
We claim that $\mathbf{u}$ and $\mathbf{y}$ satisfy

\begin{equation}
	-\mathbf{N}\mathbf{u}+\mathbf{D}\mathbf{y}=\left[
	\begin{array}
		[c]{c}%
		\mathbf{B}_\mathrm{T}\mathbf{x}_{0}\\
		0
	\end{array}
	\right]  . \label{eq11}%
\end{equation}
Substituting (\ref{eq13}) into (\ref{eq11}) we have that the first two rows of the left hand side of (\ref{eq11}) are

\begin{equation} \label{eqLHS2}
	\left[-\mathbf{N}\mathbf{D}+\mathbf{D}\mathbf{N}\right]_{(1:2m,1:k)}\bm{\theta}_{1:k} + \left[
	\begin{array}
		[c]{c}
		-\mathbf{N}_\mathrm{L}\mathbf{D}_\mathrm{U}+\mathbf{D}_\mathrm{L}\mathbf{N}_\mathrm{U}\\
		-\mathbf{N}_\mathrm{U}\mathbf{D}_\mathrm{U}+\mathbf{D}_\mathrm{U}\mathbf{N}_\mathrm{U}\\
	\end{array}
	\right]
	\mathbf{x}_0.
\end{equation}
Note that the first $2m$ rows of $\mathbf{N}\mathbf{D}$ are 
\begin{equation*}
	\left[
	\begin{array}
		[c]{cccc}%
		\mathbf{N}_\mathrm{L}\mathbf{D}_\mathrm{L} & \mathbf{0} & \mathbf{0} & \cdots\\
		\mathbf{N}_\mathrm{U}\mathbf{D}_\mathrm{L} + \mathbf{N}_\mathrm{L}\mathbf{D}_\mathrm{U} & \mathbf{N}_\mathrm{L}\mathbf{D}_\mathrm{L} & \mathbf{0} & \cdots\\
	\end{array}
	\right],  %
\end{equation*}
with a similar expression for the first $2m$ rows of $\mathbf{D}\mathbf{N}$.
By the Gohberg-Semencul formula (\ref{def:BTmatrix}), $\mathbf{D}_\mathrm{L}\mathbf{N}_\mathrm{U}-\mathbf{N}_\mathrm{L}\mathbf{D}_\mathrm{U}=\mathbf{N}_\mathrm{U}\mathbf{D}_\mathrm{L}-\mathbf{D}_\mathrm{U}\mathbf{N}_\mathrm{L}$, the matrix  multiplying $\bm{\theta}_{1:k}$ in (\ref{eqLHS2}) is identically zero. If there are more than $2m$ rows in (\ref{eq11}) then a similar argument to that used for the second block of $m$ rows can be used to show that for any such rows $k > 2m$ the left hand side of  (\ref{eq11}) is zero. Substituting $\mathbf{B}_\mathrm{T}=\mathbf{D}_\mathrm{L}\mathbf{N}_\mathrm{U}-\mathbf{N}_\mathrm{L}\mathbf{D}_\mathrm{U}$, and using $\mathbf{N}_\mathrm{U}\mathbf{D}_\mathrm{U}=\mathbf{D}_\mathrm{U}\mathbf{N}_\mathrm{U}$, we find that the second term in (\ref{eqLHS2}) is equal to the second term in (\ref{eq11}). We have shown that $\mathbf{u}$ and $\mathbf{y}$ do indeed satisfy (\ref{eq11}).  

We next show that $\mathbf{y}$ satisfies the inequality constraints for $\mathcal{P}_{z_{1:k}}(\mathbf{f}^+),$ the reasoning for validity of the inequality constraints on $\mathbf{u}$ being similar.  By (\ref{eq13}), and the first $2k$ rows of the constraints for $\bar{\mathcal{P}}(\mathbf{f}^+)$, we have
\begin{equation*}
	\mathbf{y}-\left[
	\begin{array}
		[c]{c}%
		\mathbf{N}_\mathrm{U}\mathbf{x}_{0}\\
		0
	\end{array}
	\right]
	\le \left[  \boldsymbol{1}_{4k}+\bm{\delta
	}+\bm{\gamma}\right]_{1:k} ^\mathrm{T}
\end{equation*}
and
\begin{equation*}
	-\mathbf{y}+\left[
	\begin{array}
		[c]{c}%
		\mathbf{N}_\mathrm{U}\mathbf{x}_{0}\\
		0
	\end{array}
	\right]
	\le \left[  \boldsymbol{1}_{4k}+\bm{\delta
	}+\bm{\gamma}\right]_{k+1:2k} ^\mathrm{T},
\end{equation*}
and combining these two inequalities gives  $\left\vert y_j-z_j\right\vert \leq1$ for $j=1,\ldots, k$, as required. It has been shown that $\left( \mathbf{y},\mathbf{u}\right)$ is feasible for $\mathcal{P}_{z_{1:k}}(\mathbf{f}^+)$.

We next show that $\left\langle \mathbf{x}_{k+1}\left(  \mathbf{y},\mathbf{u}%
\right)  ,\mathbf{f}^+ \right\rangle = \left\langle \theta_{k-m+1:k},\mathbf{f}^+ \right\rangle.$ This is true because \vspace{2ex} \newline
$\mathbf{x}_{k+1}\left(  \mathbf{y},\mathbf{u}\right)  =\left(  \mathbf{B}_\mathrm{T}\right)
^{-1}\left[  \mathbf{N}_\mathrm{U}u_{k-m+1:k}-\mathbf{D}_\mathrm{U}y_{k-m+1:k}\right]$
\begin{align*}
	& 	=\left(  \mathbf{B}_\mathrm{T}\right)  ^{-1}\left[  \mathbf{N}_\mathrm{U}\mathbf{D}_\mathrm{L}\theta_{k-m+1:k}-\mathbf{D}_\mathrm{U}%
	\mathbf{N}_\mathrm{L}\theta_{k-m+1:k}\right]  \text{ by (\ref{eq13})}\\
	&  =\theta_{k-m+1:k}\text{ by (\ref{def:BTmatrix}).}%
\end{align*}

For the second part of the proof we must show that any feasible solution $\left(  \mathbf{y},\mathbf{u}\right)$ to $\mathcal{P}_{z_{1:k}}(\mathbf{f}^+)$ corresponds to a feasible solution $\bm{\theta}$ to $\bar{\mathcal{P}}(\mathbf{f}^+)$, again with the same cost. Suppose $\left(  \mathbf{y},\mathbf{u}\right)$ is feasible for $\mathcal{P}_{z_{1:k}}(\mathbf{f}^+)$, implying 
\begin{equation}
	\mathbf{y}=\mathbf{D}^{-1}
	\left[\mathbf{N}\mathbf{u}+\left[
	\begin{array}
		[c]{c}%
		\mathbf{B}_\mathrm{T}\mathbf{x}_0\\
		\mathbf{0}
	\end{array} 
	\right]\right], \label{eqy}
\end{equation} 
$\left\vert y_j-z_j\right\vert \leq1$ and  $\left\vert u_j\right\vert \leq1$ for $j=1,\ldots, k$.

Put 
\begin{equation}
	\bm{\theta}=\mathbf{D}^{-1}
	\left[\mathbf{u}-\left[
	\begin{array}
		[c]{c}%
		\mathbf{D}_\mathrm{U}\mathbf{x}_{0}\\
		\mathbf{0}
	\end{array}
	\right]\right]. \label{lambda}
\end{equation}

We now show that $\bm{\theta}$ so defined satisfies the first $k$ inequalities of the constraints to $\bar{\mathcal{P}}(\mathbf{f}^+)$, namely

\begin{equation}%
	\mathbf{N} \bm{\theta} \leq  \boldsymbol{1}_{k}-\left[
	\begin{array}
		[c]{c}%
		\mathbf{N}_\mathrm{U}\mathbf{x}_{0}\\
		\mathbf{0}
	\end{array}
	\right] +\mathbf{z}.\label{Nlambda}
\end{equation}

Satisfaction of the other $3k$ inequalities can be shown using very similar arguments. Now $ \mathbf{y}-\mathbf{z}  \le \boldsymbol{1}_{k}$ implies

\begin{equation*}
	\mathbf{z} \ge \mathbf{D}^{-1}
	\left[\mathbf{N}\mathbf{u}+\left[
	\begin{array}
		[c]{c}%
		\mathbf{B}_\mathrm{T}\mathbf{x}_0\\
		\mathbf{0}
	\end{array}
	\right]\right] - \boldsymbol{1}_k,
\end{equation*}
so to demonstrate (\ref{Nlambda}) it suffices to show that

\begin{equation*}
	\mathbf{N}\mathbf{D}^{-1}
	\left[\mathbf{u}-\left[
	\begin{array}
		[c]{c}%
		\mathbf{D}_\mathrm{U}\mathbf{x}_{0}\\
		\mathbf{0}
	\end{array}
	\right]\right]
\end{equation*}
\begin{equation*}
	\le \left[
	\begin{array}
		[c]{c}%
		-\mathbf{N}_\mathrm{U}\mathbf{x}_{0}\\
		\mathbf{0}
	\end{array}
	\right] + \mathbf{D}^{-1}
	\left[\mathbf{N}\mathbf{u}+\left[
	\begin{array}
		[c]{c}%
		\mathbf{B}_\mathrm{T}\mathbf{x}_0\\
		\mathbf{0}
	\end{array}
	\right]\right], 
\end{equation*}
and in fact this is satisfied as an equality.
This follows from  the two identities $\mathbf{N}\mathbf{D}^{-1}
\equiv \mathbf{D}^{-1}\mathbf{N}$ and 
$\mathbf{D}\left[
\begin{array}
[c]{c}%
\mathbf{N}_\mathrm{U}\\
\mathbf{0}
\end{array}
\right] -
\mathbf{N}\left[
\begin{array}
[c]{c}%
\mathbf{D}_\mathrm{U}\\
\mathbf{0}
\end{array}
\right]
\equiv
\left[
\begin{array}
[c]{c}%
\mathbf{B}_\mathrm{T}\\
\mathbf{0}
\end{array}
\right].
$
The first of these identities holds because $\mathbf{N}$ and $\mathbf{D}^{-1}$ are lower triangular and toeplitz, and the first $m$ rows of the second is just a restatement of the definition of $\mathbf{B}_\mathrm{T}$.
Hence $\bm{\theta}$ is feasible for $\bar{\mathcal{P}}(\mathbf{f}^+)$.

The final step is to show equality of the cost functions, that is $ \left\langle \theta_{k-m+1:k},\mathbf{f}^+ \right\rangle = \left\langle \mathbf{x}_{k+1}\left(  \mathbf{y},\mathbf{u}%
\right)  ,\mathbf{f}^+ \right\rangle.$
We show this by demonstrating that
$\theta_{k-m+1:k} = \mathbf{x}_{k+1}\left( \mathbf{y},\mathbf{u}\right)$,
that is
\begin{equation*}
	\left[\mathbf{D}^{-1}\mathbf{u}\right]_{k-m+1:k} - \left[\mathbf{D}^{-1}
	\left[
	\begin{array}
		[c]{c}%
		\mathbf{D}_\mathrm{U}\\
		\mathbf{0}
	\end{array}
	\right]\mathbf{x}_{0}
	\right]_{k-m+1:k} 
\end{equation*}
\begin{equation}
	=\left
	(  \mathbf{B}_\mathrm{T}\right)
	^{-1}\left[  \mathbf{N}_\mathrm{U}u_{k-m+1:k}-\mathbf{D}_\mathrm{U}y_{k-m+1:k}\right], \label{eq30}
\end{equation}
where $\mathbf{y}$ is given by (\ref{eqy}).

In order to demonstrate validity of (\ref{eq30}) we need to show two things:
\begin{equation*}
	\mathbf{B}_\mathrm{T}\left[\mathbf{D}^{-1}\mathbf{u}\right]_{k-m+1:k}
\end{equation*}
\begin{equation} 
	= \mathbf{N}_\mathrm{U}u_{k-m+1:k} - 
	\mathbf{D}_\mathrm{U} \left[\mathbf{D}^{-1}\mathbf{N}\mathbf{u}\right]_{k-m+1:k} \label{eq28}
\end{equation} 
and

\begin{equation*}
	\mathbf{B}_\mathrm{T}\left[\mathbf{D}^{-1}
	\left[
	\begin{array}
		[c]{c}%
		\mathbf{D}_\mathrm{U}\\
		\mathbf{0}
	\end{array}
	\right]\mathbf{x}_{0}
	\right]_{k-m+1:k}
\end{equation*}
\begin{equation} 
	= \mathbf{D}_\mathrm{U}\left[\mathbf{D}^{-1}
	\left[
	\begin{array}
		[c]{c}%
		\mathbf{B}_\mathrm{T}\mathbf{x}_0\\
		\mathbf{0}
	\end{array}
	\right]
	\right]_{k-m+1:k}. \label{eq29}
\end{equation}

Now (\ref{eq28}) can be rewritten as
\begin{equation}
	\mathbf{B}_\mathrm{T}\left[\mathbf{D}^{-1}\mathbf{u}\right]_{k-m+1:k} + 
	\mathbf{D}_\mathrm{U} \left[\mathbf{N}\mathbf{D}^{-1}\mathbf{u}\right]_{k-m+1:k} = \mathbf{N}_\mathrm{U}u_{k-m+1:k}  
\end{equation} 
which is true because the left hand side is \newline
$
\left[\left[\mathbf{N}_\mathrm{U}\mathbf{D}_\mathrm{L}-\mathbf{D}_\mathrm{U}\mathbf{N}_\mathrm{L}\right]\mathbf{D}_\mathrm{L}^{-1} +
\mathbf{D}_\mathrm{U}\mathbf{N}_\mathrm{L}\mathbf{D}_\mathrm{L}^{-1}\right]\mathbf{u}_{k-m+1:k} 
$, and this collapses to the right hand side. 

It is straightforward to show that (\ref{eq29}) holds for any $\mathbf{x}_0$ if and only if
$
\mathbf{B}_\mathrm{T}
\mathbf{D}_\mathrm{L}^{-1}
\mathbf{D}_\mathrm{U} =
\mathbf{D}_\mathrm{U}\mathbf{D}_\mathrm{L}^{-1}\mathbf{B}_\mathrm{T} 
$, and this is an identity by virtue of (\ref{def:BTmatrix}).

This completes the proof of the equivalence of $\bar{\mathcal{P}}(\mathbf{f}^+)$ and 
$\mathcal{P}_{z_{1:k}}(\mathbf{f}^+),$ and it follows that $\mathcal{P}(\mathbf{f}^+)$ and $\mathcal{D}_{z_{1:k}}(\mathbf{f}^+)$ are a dual pairing in the sense of (\ref{asymdual}). The final step is to show that the complementary slackness conditions, (i) and (ii), imply the alignment conditions of the theorem statement. This can be done using (\ref{def:BTmatrix}), (\ref{notation_dual}), (\ref{eqy}) and (\ref{lambda}). The aglebraic manipulations involved are similar to those already used above, and are omitted.

\section{Proof of Theorem \ref{main1}}
We are given $\mathbf{x}\in S_{k}$ and
$\mathbf{f}\in \mathcal{N}_{S_{k}}(\mathbf{x})$. Since $S_{k}$ is non-empty, there exists $\big(y_{1:k-1},u_{1:k-1} \big) \in \arg\max \mathcal{P}_{z_{1:k-1}}(\mathbf{f})$ and $\big( y^{\ast}_{1:k-1},u^{\ast}_{1:k-1} \big) \in \arg\min \mathcal{D}_{z_{1:k-1}}(\mathbf{f})$, and these sequence pairs are aligned by Theorem \ref{main}. 

We first prove the if part of the theorem statement. Given $\left(
u_{k},y_{k}^{\ast}\right)  \in M\left( \mathbf{x},  \mathbf{f},z_{k}\right)$ we put $y_{k}  =\mathbf{C}\mathbf{x}+Du_k$ and $u_{k}^{\ast} =\mathbf{C}^{\ast}\mathbf{f}+D^{\ast}y_{k}^{\ast}$. By Definition \ref{defM}, $\left(y_k,u_k\right)$ is aligned with $\left(
y_k^{\ast},u_k^{\ast}\right)$. Put $\mathbf{f^+}
=\mathbf{A}^{\ast}\mathbf{f}+\mathbf{B}^{\ast}y_{k}^{\ast
}$. Then
$\big(y_{1:k},u_{1:k} \big)$ and $\big(y_{1:k}^*,u_{1:k}^* \big)$ are feasible, respectively, for $\mathcal{P}_{z_{1:k}}(\mathbf{f}^+)$ and $\mathcal{D}_{z_{1:k}}(\mathbf{f}^+)$, and are aligned. A further application of Theorem \ref{main} gives $\big(y_{1:k},u_{1:k} \big) \in \arg\max \mathcal{P}_{z_{1:k}}(\mathbf{f}^+)$ and $\big( y^{\ast}_{1:k},u^{\ast}_{1:k} \big) \in \arg\min \mathcal{D}_{z_{1:k}}(\mathbf{f}^+)$. Then, by Proposition \ref{propconeargmax} (with $k$ replaced by $k+1$), $\mathbf{f}^+\in \mathcal{N}_{S_{k+1}}\left(\mathbf{x}^+ \right)$, where $\mathbf{x}^+ =\mathbf{A}\mathbf{x}+\mathbf{B}u_{k}$.

For the only if part, $\mathbf{x}^+ =\mathbf{A}\mathbf{x}+\mathbf{B}u_{k} \in S_{k+1}$ and  $\mathbf{f^+}
=\mathbf{A}^{\ast}\mathbf{f}+\mathbf{B}^{\ast}y_{k}^{\ast}\in \mathcal{N}_{S_{k+1}}(\mathbf{x^+})$ are given. Put $y_{k}  =\mathbf{C}\mathbf{x}+Du_k$ and $u_{k}^{\ast} =\mathbf{C}^{\ast}\mathbf{f}+D^{\ast}y_{k}^{\ast}$. There exist sequences 
$\big(y_{1:k},u_{1:k} \big) \in \arg\max \mathcal{P}_{z_{1:k}}(\mathbf{f}^+)$ and $\big( y^{\ast}_{1:k},u^{\ast}_{1:k} \big) \in \arg\min \mathcal{D}_{z_{1:k}}(\mathbf{f}^+)$ and by Theorem \ref{main} they are aligned. Hence $\left(y_k,u_k\right)$ is aligned with $\left(y_k^{\ast},u_k^{\ast}\right)$. The inequalities $\left\vert u_k\right\vert \le 1$ and $\left\vert y_k-z_{k}\right\vert \le 1$, as well as the condition  $(y_k,u_k) \in L\left(  \mathbf{x}\right)$, hold because the constraints to $\mathcal{P}_{z_{1:k}}(\mathbf{f}^+)$ are satisfied by $\big(y_{1:k},u_{1:k} \big)$. Also $(y_k^{\ast},u_k^{\ast}) \in L^{\ast}\left(  \mathbf{f}\right)$ because the constraints to $\mathcal{D}_{z_{1:k}}(\mathbf{f}^+)$ are satisfied by $(y_{1:k}^{\ast},u_{1:k}^{\ast})$. All of the conditions of Definition \ref{defM} are satisfied, so $\left(
u_{k},y_{k}^{\ast}\right)  \in M\left( \mathbf{x}, \mathbf{f},z_{k}\right)$, as required.

\section{Proof of Theorem \ref{main2}}
By the definition of precursor,  $\mathbf{x} = \mathbf{x}_{k}\left(  y_{1:k-1},u_{1:k-1}\right) $ for some $\left( y_{1:k},u_{1:k}\right) $ feasible for $\mathcal{P}_{z_{1:k}}(\cdot)$. From the state space description of the primal system $\mathbf{x}^+=\mathbf{A}\mathbf{x}+Bu_k$ and, by Proposition \ref{propconeargmax}, $\left( y_{1:k},u_{1:k}\right) \in \arg\max\mathcal{P}_{z_{1:k}}\left(\mathbf{f}^+ \right)$ for all $\mathbf{f}^+\in \mathcal{N}_{S_{k+1}}(\mathbf{x}^+)$. Note that $\mathcal{N}_{S_{k+1}}(\mathbf{x}^+)$ is never empty because $\mathcal{N}_{S_{k+1}}(\mathbf{x}^+)$ must contain the zero vector. Now for all $\mathbf{f}^+\in \mathcal{N}_{S_{k+1}}(\mathbf{x}^+)$ there exists $\left(  y_{1:k}^{\ast},u_{1:k}^{\ast}\right) \in \arg\min\mathcal{D}_{z_{1:k}%
}\left(  \mathbf{f}^+\right)$ and, by the state space description of the dual system, $\mathbf{f}^+=A^*\mathbf{f}+B^*y_k^*$, where $\mathbf{f} =  \mathbf{x}_{k}^*\left(y^*_{1:k-1},u^*_{1:k-1}\right)$. By Theorem \ref{main}, $\left( y_{1:k},u_{1:k}\right)$ and $\left( y_{1:k}^*,u_{1:k}^*\right)$ are aligned. We have shown that $\left( y_{1:k-1},u_{1:k-1}\right)$ is feasible for $\mathcal{P}_{z_{1:k-1}}(\mathbf{f})$, and is aligned with $\left( y_{1:k-1}^{\ast},u_{1:k-1}^{\ast}\right)$, which is feasible for $\mathcal{D}_{z_{1:k-1}}(\mathbf{f})$. By Theorem \ref{main}, $\left( y_{1:k-1},u_{1:k-1}\right) \in \arg\max \mathcal{P}_{z_{1:k-1}}(\mathbf{f})$ and, from Proposition \ref{propconeargmax}, we have $\mathbf{f} \in \mathcal{N}_{S_{k}}\bigl(\mathbf{x}_{k}\left(y_{1:k-1},u_{1:k-1}\right) \bigr) =  \mathcal{N}_{S_{k}}\bigl(\mathbf{x} \bigr)$.

Finally we show that $\left(
u_{k},y_{k}^{\ast}\right)  \in M\left( \mathbf{x},  \mathbf{f},z_{k}\right)$ for any $\left(
u_{k},y_{k}^{\ast}\right)$ as given in the previous paragraph. From the state space descriptions of the primal and dual systems, respectively, we have $(y_k,u_k) \in L\left(\mathbf{x}\right)$ and $(y_k^*,u_k^*) \in L^*\left(\mathbf{f}\right)$. Furthermore 
$(y_k,u_k)$ and $(y_k^*,u_k^*)$ are aligned, because by Theorem \ref{main}, $\left( y_{1:k},u_{1:k}\right)$ and $\left( y_{1:k}^*,u_{1:k}^*\right)$ are aligned. All the conditions of Definition \ref{defM} are satisfied, so
$\left(
u_{k},y_{k}^{\ast}\right)  \in M\left( \mathbf{x},  \mathbf{f},z_{k}\right)$.

\end{document}